\documentclass[10pt, a4paper]{article} 
\author{
Rafael Granero-Belinch\'on
\\{\footnotesize Departamento de Matem\'aticas, Estad\'istica y Computaci\'on}
\\{\footnotesize Universidad de Cantabria}
\\{\footnotesize Santander, Espa\~na}
\\{\footnotesize email: {\it rafael.granero@unican.es}}
\and 
Stefano Scrobogna
\\{\footnotesize Departamento de An\'alisis Matemático}
\\{\footnotesize Universidad de Sevilla}
\\{\footnotesize Sevilla, Espa\~na}
\\{\footnotesize email: {\it scrobogna@us.es}}
}
\usepackage{amsfonts}

\usepackage[english]{babel}
\usepackage{pst-grad} % For gradients
\usepackage{pst-plot} % For axes
\usepackage{pstricks}
\usepackage{amsmath,amssymb,mathrsfs,amsthm, mathtools}
\usepackage{stmaryrd}
\usepackage{tikz} 
\usepackage{mathcomp,wasysym}  
\usepackage{graphicx}  
\usepackage[all]{xy} \xyoption{arc} \xyoption{color}
\usepackage{epsfig}
\usepackage{cite}
\usepackage{enumerate}
\usepackage[a4paper,
left=2.5cm, right=2.5cm,
top=3cm, bottom=3cm]{geometry}
\usepackage[a4paper]{geometry}
\usepackage{empheq}
\usepackage{times}
\usepackage{accents}

\DeclareMathAlphabet{\mathcal}{OMS}{cmsy}{m}{n}

\DeclareFontFamily{U}{mathc}{}
\DeclareFontShape{U}{mathc}{m}{it}%
{<->s*[1.03] mathc10}{}

\DeclareMathAlphabet{\mathscr}{U}{mathc}{m}{it}

\newcommand{\dx}{\textnormal{d}{x}}
\newcommand{\dy}{\textnormal{d}{y}}
\newcommand{\dt}{\textnormal{d}{t}}
\newcommand{\ddt}{\frac{\textnormal{d}}{ \textnormal{d}{t}}}

\newcommand{\dd}{\textnormal{d}}

\newcommand{\sgn}{\textnormal{sgn}}
\newcommand{\pare}[1]{\left( #1 \right)}
\newcommand{\norm}[1]{\left\| #1 \right\|}
\newcommand{\av}[1]{\left| #1 \right|}
\newcommand{\bra}[1]{\left[ #1 \right]}

\newcommand{\set}[1]{\left\{ #1 \right\}}

\newcommand{\cA}{\mathcal{A}}
\newcommand{\sA}{\mathsf{A}}
\newcommand{\cC}{\mathcal{C}}

\newcommand{\cS}{\mathcal{S}}

\newcommand{\bR}{\mathbb{R}}
\newcommand{\bT}{\mathbb{T}}
\newcommand{\bZ}{\mathbb{Z}}
\newcommand{\bN}{\mathbb{N}}

\newcommand{\ck}{\mathscr{k}}

\newcommand{\pat}{\partial_t}

\newcommand{\jeps}{\mathcal{H}_\kappa*}

\normalsize
\normalsize
\usepackage[bookmarks,colorlinks]{hyperref}

\newcommand{\loc}{\textnormal{loc}}

\theoremstyle{theorem}
\newtheorem{theorem}{Theorem}[section]
\newtheorem*{theorem*}{Theorem}
\newtheorem{prop}[theorem]{Proposition}
\newtheorem{lemma}[theorem]{Lemma}

\theoremstyle{definition}

\numberwithin{equation}{section}

\title{Well-posedness of the water-wave with viscosity problem}

\begin{document}

%\subjclass{35R35, 35Q35, 35S10, 76B03}
%\keywords{Water waves, damping, moving interfaces, free-boundary problems}

%\ccode{2010 AMS Subject Classification: 35455, 35B41, 92C17}

\maketitle

\begin{abstract}
In this paper we study the motion of a surface gravity wave with viscosity. In particular we prove two well-posedness results. On the one hand, we establish the local solvability in Sobolev spaces for arbitrary dissipation. On the other hand, we establish the global well-posedness in Wiener spaces for a sufficiently large viscosity. These results are the first rigorous proofs of well-posedness for the Dias, Dyachenko \& Zakharov system ({\em Physics Letters A} 2008) modelling gravity waves with viscosity when surface tension is not taken into account.
	\end{abstract}
{\small
\tableofcontents}

\allowdisplaybreaks
\section{Introduction}
The motion of water waves is a classical research topic that has attracted a lot of attention from many different researchers in Mathematics, Physics and Engineering. For most applications, a gravity surface wave in deep water is described by the following free boundary problem
\begin{align}
\label{euler}
\rho\left(\pat u+(u\cdot\nabla)u\right)+\nabla p+G\rho e_2&=0\quad x\in \Omega(t),\,t\in[0,T]\\
\partial_t \rho+ u\cdot\nabla \rho&=0\quad x\in \Omega(t),\,t\in[0,T]\\
\nabla\cdot u&=0\quad x\in \Omega(t),\,t\in[0,T]\\
\nabla\times u&=0\quad x\in \Omega(t),\,t\in[0,T],
\end{align}
where, $\Omega(t)\subset \mathbb{R}^2$ is the moving spatial domain, $T>0$ and $u=(u_1,u_2)$, $\rho$, $p$ and $G$ denote the incompressible velocity field, the density of the fluid, the pressure and the acceleration due to gravity, respectively. The previous system can also be equivalently written in terms of the velocity potential $u=\nabla \Phi$:
\begin{equation}\label{euler2}
\left\lbrace
\begin{aligned}
& \Delta \Phi=0&&\text{ in }\Omega(t)\times[0,T],\\
%%--------------------------------------------------------
%&\phi =\xi &&\text{ on }\Gamma(t)\times[0,T],\\
%--------------------------------------------------------
& \lim_{x_2\rightarrow-\infty}\nabla\Phi= 0 && \forall\,x_1, \\
%--------------------------------------------------------
& \rho\left(\Phi_t+\frac{1}{2}|\nabla\Phi|^2+ Gh\right) =0&&\text{ on }\Gamma(t)\times[0,T],\\
%--------------------------------------------------------
& h_t =\nabla\Phi\cdot(-\partial_1 h,1)&&\text{ on }\Gamma(t)\times[0,T],
\end{aligned}
\right. 
\end{equation}
Here water is assumed to be incompressible, irrotational and inviscid. Usually, these assumptions are enough to describe the main part of the dynamics of real water waves. However, sometimes viscosity needs to be taken into account. For instance, Bona, Pritchard \& Scott \cite{bona1981evaluation}  wrote that 

\begin{quote}
For these five experiments it was found that the inclusion of a dissipative term was much more im portant than the inclusion of the nonlinear term,  although the inclusion of the nonlinear term was undoubtedly beneficial in describing the observations.
\end{quote}

We refer to \cite{GS19} for more examples of the necessity of the viscosity effects in experiments.

When viscosity is considered, the free boundary Euler equations \eqref{euler} need to be replaced by the free boundary Navier-Stokes equations and the irrotationality assumption should be removed. This will lead to a fairly complicated free boundary problem. However, under certain hypotheses, it is well-known that vorticity only plays a role in a layer near the free boundary. The observation of this fact dates back at least to the classical works of Lamb \cite{lamb1932hydrodynamics} and Boussinesq \cite{boussinesq1895lois} in the XIX century. Thus, it would be desirable to add dissipative effects to \eqref{euler} (equivalently \eqref{euler2}) without going all the way to the Navier-Stokes equations and the subsequent removal of the irrotationality assumption. 

In a recent paper, Dias, Dyachenko \& Zakharov \cite{DDZ08} proposed a system to solve that question. These authors started with the linearized Navier-Stokes equations and derived the following system for the motion of a surface wave under gravity and viscous effects 
\begin{equation}\label{DDZ}\tag{DDZ}
\left\lbrace
\begin{aligned}
& \Delta \Phi=0&&\text{ in }\Omega(t)\times[0,T],\\
%%--------------------------------------------------------
%&\phi =\xi &&\text{ on }\Gamma(t)\times[0,T],\\
%--------------------------------------------------------
& \lim_{x_2\rightarrow-\infty}\nabla\Phi= 0 && \forall\,x_1, \\
%--------------------------------------------------------
& \rho\left(\Phi_t+\frac{1}{2}|\nabla\Phi|^2+ Gh\right) =-2\nu \partial_2^2\Phi&&\text{ on }\Gamma(t)\times[0,T],\\
%--------------------------------------------------------
& h_t =\nabla\Phi\cdot(-\partial_1 h,1)+2\frac{\nu}{\rho}\partial_{1}^2h &&\text{ on }\Gamma(t)\times[0,T],
\end{aligned}
\right. 
\end{equation}
where 
\begin{align}\label{Omegat}
\Omega(t) & = \set{ (x_1, x_2)\in\bR^2 \ \Big| \ {-L}\pi <x_1<{L}\pi\,, -\infty < x_2 < h(x_1,t)\,, \ t\in[0,T] }, \\
%----------------------------------
\label{Gammat}
\Gamma(t) & = \set{ \pare{x_1,h(x_1,t)}\in\bR^2 \ \Big| \ {-L}\pi <x_1<{L}\pi\,,\ t\in[0,T] },
\end{align} 
and $2\pi L$ is the characteristic wavelength of the wave. In the previous system (that we call DDZ system) $\nu$ denote the the dynamic viscosity of the fluid. In this system, the free surface experiments dissipative effects and the fluid is assumed to be irrotational. In other words, the velocity field is given by
$$
u=\nabla\Phi.
$$
The previous system has to be supplemented with zero mean initial data
$$
h(x_1,0)=h_0(x_1),
$$
$$
\Phi(x_1,h(x_1,0),0)=\xi_0(x_1).
$$

It is useful to write \eqref{DDZ} in dimensionless form. Following \cite{GS19}, we denote by $H$ and $L$ the typical amplitude and wavelength of the water wave. We change to dimensionless variables (denoted with $\tilde{\cdot}$)
\begin{align*}
x=L \ \tilde{x}, & \quad\quad t=\sqrt{\frac{L}{G}}\ \tilde{t}, \\
h(x_1,t)=H  \ \tilde{h}(\tilde{x}_1,\tilde{t}), & \quad\quad \Phi(x_1,x_2,t)=H\sqrt{G L}\tilde{\Phi}(\tilde{x}_1,\tilde{x}_2,\tilde{t}).
\end{align*}
with the dimensionless fluid domain
\begin{align*}
\widetilde{\Omega}(t) & =\set{ \pare{ \tilde{x}_1, \tilde{x}_2} \ \left| \ -\pi<\tilde{x}_1< \pi\,, -\infty < \tilde{x}_2 < \frac{H}{L}\tilde{h}(\tilde{x}_1,t)\,, \ t\in[0,T]\right.  }, \\
%--------------------------------
\widetilde{\Gamma}(t) & =\set{ \pare{ \tilde{x}_1, \frac{H}{L}\tilde{h}(\tilde{x}_1,t)}\,, \ t\in[0,T] }
\end{align*}
We define the following dimensionless parameters:
\begin{align} \label{eq:dimensionless_parameters}
\varepsilon=\frac{H}{L}, && \alpha=\frac{2\nu}{\rho \sqrt{G}L^{\frac{3}{2}}}.
\end{align}
Dropping the tildes in order to simplify the notation, and using the notation
$$
\xi(x_1,t)=\Phi(x_1,h(x_1,t),t),
$$
we find the dimensionless form of the damped water waves problem \eqref{DDZ}:
\begin{subequations}\label{eq:all2dimensionless}
\begin{align}
\Delta \Phi&=0&&\text{ in }\Omega(t)\times[0,T],\\
\Phi  &= \xi \qquad &&\text{on }\Gamma(t)\times[0,T],\\
\lim_{x_2\rightarrow-\infty}\nabla \Phi&= 0 \qquad &&\text{on }\bT\times[0,T],\\
\xi_t &=-\frac{\varepsilon}{2}|\nabla \Phi|^2- h-\alpha\partial_2^2\Phi &&\nonumber\\
&\quad+\varepsilon\partial_2\Phi\left(\nabla\Phi\cdot(-\varepsilon\partial_1 h,1)+\alpha \partial_1^2 h\right)&&\text{ on }\Gamma(t)\times[0,T],\\
h_t&=\nabla\Phi\cdot (-\varepsilon\partial_1 h,1)+\alpha\partial_{1}^2 h&&\text{ on }\Gamma(t)\times[0,T],\\
\xi &=\xi_0&&\text{ on }\bT,\\
h&=h_0&&\text{ on }\bT.
\end{align}
\end{subequations}
This formulation of the water waves with viscosity problem recovers the Zakharov/Craig-Sulem \cite{zakharov1968stability,craig1993numerical} formulation of the water wave problem when $\alpha=0$. The closely related problem of finding asymptotic models that describe the motion of water waves with viscosity has been considered by many different researchers. For instance, Dutykh \& Dias \cite{dutykh2007viscous} obtain a nonlocal (in time) Boussinesq system (see also \cite{dutykh2009visco, dutykh2007dissipative} and the references therein), Kakleas \& Nicholls \cite{kakleas2010numerical} derived the viscous analog of the classical Craig-Sulem WW2 model (see also \cite{ambrose2012well}) and the authors derived a single nonlocal wave equation in \cite{GS19} (see also \cite{gs2019well}).

To the best of the authors knowledge, the only mathematical result addressing the well-posedness of the water waves with viscosity free boundary problem (in the spirit of the Dias, Dyachenko \& Zakharov model) is the recent paper by Ngom \& Nicholls \cite{ngom2018well}. In such work the authors proved that given an initial data $ \epsilon\pare{h_0\pare{x_1}, \xi_0\pare{x_1}}, \ 0 < \epsilon \ll 1 $ and a non-vanishing surface tension $\sigma\neq 0$ the system \eqref{eq:all2dimensionless} admit a unique global solution which moreover becomes instantaneously analytic in the amplitude parameter in the sense that
\begin{equation*}
\pare{h\pare{x_1, t}, \xi\pare{x_1, t}} = \epsilon \sum_{n=0}^\infty \epsilon^n \pare{h_n\pare{x_1, t}, \xi_n\pare{x_1, t}}, 
\end{equation*} 
for any $ t > 0 $. 
We would like to remark that in the present manuscript we restrict ourselves to the case of vanishing surface tension $\sigma=0$ and, as a consequence, our results are outside the framework of \cite{ngom2018well}. In addition in Theorem \ref{theorem2} we prove as well that the solution becomes instantaneously analytic in the \textit{tangential direction} to the interface, which is again another difference when compared to \cite{ngom2018well}. Furthermore, the regularity required for the initial data is almost sharp in the sense that the initial data belongs to a functional space which enjoys the same scaling invariance as $ W^{1, \infty} $ (see Appendix \ref{sec:notation} for the detailed description of the functional spaces considered in the manuscript). This regularity threshold is required for the interface in order the strong solution to elliptic problem can be found. 

Let us now introduce a detailed statement of the results proved in the present work, first we prove the local existence of classical solution in Sobolev class for arbitrary $\alpha,\varepsilon>0$ (see Appendix \ref{sec:notation} for the definition of the functional spaces that we use in this work):
\begin{theorem}\label{theorem1}There exist small enough constants $0<\gamma,T_*<\infty$ such that if the initial data $(h_0,\xi_0)\in H^3(\bT)\times H^3(\bT)$ satisfy
$$
\max\{|h_0|_2,|\xi_0|_2\}<\gamma
$$
then there exists a smooth solution
$$
h\in C([0,T_*],H^3(\bT))\cap L^2(0,T;H^4(\bT)),
$$
$$
\xi\in C([0,T_*],H^3(\bT))\cap L^2(0,T;H^4(\bT)),
$$
to the damped water wave problem \eqref{DDZ}.
\end{theorem}

Although it seems that the DDZ system reduces to a viscous regularization of the hyperbolic nonlocal problem \eqref{euler2}, the fact that we are dealing with a free boundary problem makes the previous statement (at least) unaccurate. Actually \eqref{DDZ} is a cross diffusion system and this fact makes necessary the size restriction on the initial data. Let us try to motivate the requirement of the size constraint in the initial data. In order to do that we consider the following toy problem
\begin{equation*}
\begin{small}
\begin{aligned}
\xi_t &=\alpha\Lambda\xi \partial_1^2 h-h+\alpha \partial_1^2 \xi,\\
h_t&=-\partial_1 h \partial_1 \xi+\Lambda \xi +\alpha  \partial_1^2 h,
\end{aligned}
\end{small}
\end{equation*}
where $\Lambda=(-\partial_1^2)^{1/2}$. Then we see that in order the contributions from the term $\alpha\Lambda\xi \partial_1^2 h$ can be absorbed by the parabolic term $\alpha  \partial_1^2 h$ we need a smallness condition in $\Lambda \xi$.

Let us explain the methodology behind the proof of Theorem \ref{theorem1}. First we fixed the domain using a diffeomorphism with optimal Sobolev regularity. This is known in the literature as using the Arbitrary Lagrangian-Eulerian formulation \cite{CGS,hadvzic2015global}. This formulation allows us to control the velocity potential $\Phi$ in the bulk of the fluid while keeping track of the regularity of the interface $\xi$ and the trace of the velocity potential $\xi$. Then we regularize appropriately the problem in such a way the energy estimates are independent of the parameter of regularization. As it is well-known, this is a rather challenging issue when dealing with free boundary problems. Finally we perform $L^2$ based energy estimates in order we have the required estimates that allow us to pass to the limit in the regularization parameter.

After establishing the local well-posedness for arbitrary parameter $\alpha$, we prove the global existence and instantaneous gain of analyticity for small initial data in Wiener spaces if $\alpha$ is beyond a threshold value (see Appendix \ref{sec:notation} for the definition of the functional spaces that we use in this work):
	\begin{theorem} \label{theorem2}
Let $\alpha>2$ and $\varepsilon>0$ be two fixed constants. There exists a constant $\mathcal{C}(\alpha,\epsilon)$ such that if $h_0 \in \mathbb{A}^1_0 $, $\xi_0 \in \mathbb{A}^1_0 $ are the zero mean initial data satisfying
\begin{align}\label{eq:smllnss0}
\av{h_0}_{1} + \av{\xi_0}_{1} < \mathcal{C}\ll1, 
\end{align}
then there exists a global solution of \eqref{DDZ}. Furthermore, this solution becomes instantaneously analytic in a growing strip in the complex plane. In particular, for any
\begin{equation*}
\mu\in [0,\alpha/2),	
\end{equation*}
the solution lies in the energy space
\begin{align*}
h \in \cC\pare{\bra{0, T}; \mathbb{A}^1_{\mu t}}\cap L^1\pare{\bra{0, T}; \mathbb{A}^3_{\mu t}},\forall\,0<T<\infty. 
\end{align*}
In addition, the solution decays
$$
\pare{\av{\xi}_{1, \mu t}  + \av{h}_{1, \mu t}}\leq e^{-\delta t}\pare{\av{\xi_0}_{1}  + \av{h_0}_{1}},
$$
for certain $0<\delta=\delta(\alpha,\varepsilon)\ll1$.
\end{theorem}
We would like to emphasize that the condition
$$
\alpha>2
$$
is not optimal. We can improve it to
$$
\alpha>1
$$
but the proof will be more technical without leading to a great improvement in the result. Let us explain the reason behind the assumption on $\alpha$. To do that we consider the linear problem in Fourier variables
\begin{equation*}
\begin{small}
\begin{aligned}
\hat{\xi}_t &=-\hat{h}-\alpha k^2 \hat{\xi},\\
\hat{h}_t&=|k|\hat{\xi} -\alpha k^2 \hat{h}.
\end{aligned}
\end{small}
\end{equation*}
Then we see that the $\ell^1$ norm evolve according to (see Appendix \ref{sec:notation} for the definition of the Wiener spaces)
\begin{equation*}
\begin{small}
\begin{aligned}
\frac{d}{dt}\|\hat{\xi}\|_{\ell^1} &\leq \|\hat{h}\|_{\ell^1}-\alpha  \|k^2\hat{\xi}\|_{\ell^1},\\
\frac{d}{dt}\|\hat{h}\|_{\ell^1}&=\||k|\hat{\xi}\|_{\ell^1} -\alpha  \|k^2\hat{h}\|_{\ell^1}.
\end{aligned}
\end{small}
\end{equation*}
Then, in order the contributions from the $\alpha-$independent terms can be absorbed into the dissipative part of the linear problem, we need $\alpha>1$. Let us try to further ellaborate on the methodology of the proof of Theorem \ref{theorem2}. In order to establish this result, we use the Wiener-Sobolev spaces introduced by the authors in collaboration with Gancedo in \cite{GGS19} (we refer to \ref{sec:notation} for the precise definition). These non-standard spaces are \emph{anisotropic} in the sense that they measure differently the $x_1$ and the $x_2$ variables. In particular, they measure $x_1$ akin to a Wiener space and $x_2$ as a $W^{k,1}$ Sobolev space. However, this makes them very useful when dealing with a free boundary problem because it allow us to control the interface in a certain Wiener spaces while keeping track of the regularity in the bulk of the fluid in a Sobolev sense.

\subsection{Notation}\label{sec:notationa}
All along the manuscript $ \bT = \bR / 2\pi\bZ $, which alternatively can be thought as the interval $ \bra{-\pi, \pi} $ endowed with periodic boundary conditions. Given a matrix $ A \in \bR^{n\times m} $ we denote with $A_j^i $ the entry of $ A $ at row $ i $ and column $ j $, and we use Einstein convention for the summation of repeated indexes  from this point onwards. We write
$$
\partial_j f=\frac{\partial f}{\partial x_j},\quad f_t=\frac{\partial f}{\partial t }
$$
for the space derivative in the $j-$th direction and for a time derivative, respectively. We also write
$$
f,_j=\frac{\partial f}{\partial x_j},
$$
and
$$
D^\alpha u =\partial_1^{\alpha_1}\partial_2^{\alpha_2}u,\;\alpha_1+\alpha_2=\alpha.
$$
Let $v(x_1)$ denote a $L^2$ function on $\bT$ (as usually, identified with the interval $[-\pi,\pi]$ with periodic boundary conditions). We recall its Fourier series representation:
	\begin{equation*}
	\hat{v}\pare{n} = \frac{1}{2\pi}\int _{\bT} v\pare{x} e^{-in x}\dx, 
	\end{equation*}
where $ n\in\bZ $. Then we have that
\begin{equation*}
	v\pare{x} =\sum_{n\in\bZ} \hat{v}\pare{n} \ e^{in x}. 
	\end{equation*}
We define the Calderon operator $\Lambda $ as
	\begin{equation*}
	\widehat{\Lambda v}\pare{n} = \av{n} \hat{v}\pare{n}.
	\end{equation*}
Similarly, for sufficiently smooth functions $f$
	\begin{equation*}
	\widehat{f\pare{\Lambda} v}\pare{n} = f\pare{\av{n}} \hat{v}\pare{n}.
	\end{equation*}
We define
\begin{equation*}
\mathcal{H}_\kappa(x_1)
\end{equation*}
as the periodic heat kernel at time $t=\kappa$, and we denote
$$
f^\kappa=\jeps f \ \text{ and } \ f^{\kappa\kappa}=\jeps\left(\jeps f\right)\,.
$$

We write $C$ for a constant that may change from line to line. This constant may depend on the norm of the initial interface and the initial velocity at the free boundary in lower norms $C=C(|h_0|_{2},|\xi_0|_{2})$. We write $\mathcal{P}$ for a polynomial that may change from line to line. For the definition and some properties of the functional spaces used in this paper we refer to the Appendix \ref{sec:notation}

\section{The Arbitrary Lagrangian-Eulerian formulation of the damped water waves system}
\subsection{A family of diffeomorphisms}
Our analysis of the damped water waves problem \eqref{DDZ} is based on a time-dependent change-of-variables $\psi(x_1,x_2,t)$. This family of diffeomorphisms flattens the domain and changes the free boundary problem to one set on the smooth reference domain
$$
\Omega= \bT \times\pare{ -\infty, 0}.
$$
The reference boundary is then given by
$$
\Gamma=\bT\times\{0\}.
$$
We let $N=e_2$ denote the outward pointing unit normal on $\Gamma$. We observe that we are dealing with a parabolic problem with cross-diffusion. Thus our analysis relies on obtaining a parabolic gain of regularity. To do so, we need a reference domain $\Omega$ with smooth boundary.

In this section we construct a family of diffeomorphism with optimal Sobolev regularity. To do that we follow the ideas in \cite{CGS}. The main point of this construction is to solve appropriate Poisson problems to find $\psi(x_1,x_2,t)$. Then, an application of the inverse function theorem together with standard elliptic estimates will show that these mappings $\psi(x_1,x_2,t)$ are a family of diffeomorphisms with the desired smoothness (see \cite{GGS19} for more details). Given a function $h\in C(0,T;H^{2})$ with initial data $h(x_1, 0)=h_0(x_1)$, we define $(\delta \psi)$ as the solution of the following elliptic system:
\begin{subequations}\label{psita}
\begin{alignat*}{2}
\Delta \delta\psi&=0  \qquad&&\text{in}\quad \Omega\times[0,T]\,,\\
\delta\psi&= h \qquad &&\text{on}\quad \Gamma\times[0,T]\,.
\end{alignat*}
\end{subequations}
Due to elliptic estimates, we have the bound
\begin{align}\label{closetopsi0}
\|\nabla\delta\psi(t)\|_{r-0.5}&\leq C|h(t)|_{r}\text{ for all }r>0.5.  
\end{align}
Then we define the mapping
$$
\psi(x_1,x_2,t)=e+(0,\delta\psi)=(x_1,x_2+\delta\psi(x_1,x_2,t)).
$$
We observe that, if $h$ is small enough in certain sense, then $\psi$ is a diffeomorphism. To see that we use that its distance from the identity map, $e(x)$, is small (see \cite{CGS,GGS19}). 

As this paper consider the case where the initial domain lies in Sobolev and Wiener spaces, we need to show that $\psi$ is a diffeomorphism in both functional settings. Indeed, let us first assume that $h$ is small in certain Sobolev norms, namely
$$
|h|_{2}\ll1.
$$
Then,
$$
\|\psi-e\|_{\cC^1}\leq \|\delta\psi\|_{\cC^1}\leq C\|\nabla\delta\psi\|_{1.25}\leq C|h|_{1.75}\ll1.
$$
Thus, we have that
$$
\|\psi(x)-\psi(y)\|_{\cC^0}\geq\|x-y\|_{\cC^0}-\|\delta\psi(x)-\delta\psi(y)\|_{\cC^0}\geq (1-C|h|_{1.75})\|x-y\|_{\cC^0},
$$
which ensures that $\psi$ is injective if $h$ is small enough in $H^{2}(\bT)$. As a consequence of the inverse function theorem $ \psi $ is a (global) diffeomorphism such that
$$
\|\nabla\psi(t)-\text{Id}\|_{2.5}<\infty.
$$

\subsection{The matrix $\sA$}\label{sec:matA}
We write
$$
J=\text{det}(\nabla\psi)=\psi^1,_1\psi^2,_2-\psi^2,_1\psi^1,_2, = 1+ \delta\psi,_2.
$$	
We have the bound
\begin{align}\label{J1.25}
\|J(t)-1\|_{1.25}=\|\delta\psi,_2\|_{1.25}\leq C|h(t)|_{1.75}. 
\end{align}

We write 
$$
\sA = (\nabla \psi)^{-1} = J^{-1} \left(\begin{array}{cc}
\psi^2,_2 & - \psi^1,_2 \\
-\psi^2,_1 & \psi^1,_1
\end{array}
\right)=\frac{1}{1+\delta\psi,_2} \left(\begin{array}{cc}
1+\delta\psi,_2 & 0 \\
-\delta\psi,_1 & 1
\end{array}
\right)\,.
$$
Using the fact that
$$
\sA^i_r\psi^r,_{j}=\delta^i_j,
$$
we obtain the useful identities
\begin{equation}\label{propA}
(\sA_t)^i_k=-\sA^i_r(\psi_t)^r,_j \sA^j_k,\qquad D^2\sA=-2 D^1\sA\nabla D^1\psi \sA-\sA D^2\nabla\psi \sA.
\end{equation}
We will also make use of the Piola's identity:  
$$
(J\sA^k_i),_k=0.
$$

\subsection{The ALE problem in the fixed domain}
We can now define the pull-back of the scalar potential $ \Phi $ as 
\begin{equation*}
\phi = \Phi \circ \psi.
\end{equation*}
We compute that
$$
(\nabla \Phi)\circ \psi=\sA^T\nabla\phi=\sA^k_i\phi,_k,
$$
and
$$
(\Delta \phi)\circ\psi=\sA^i_j(\sA^k_j\phi,_k),_i.
$$
Then, defining
$$
n=\frac{(-\varepsilon h,_1,1)}{\sqrt{1+(\varepsilon h,_1)^2}},\;g=1+(\varepsilon h,_1)^2,\;\tilde{n}=\sqrt{g}n=(-\varepsilon h,_1,1),
$$
the problem \eqref{eq:all2dimensionless} can be written as
\begin{subequations}\label{eq:ALE}
\begin{align}
\sA^\ell_j(\sA^k_j\phi,_k),_\ell&=0&&\text{ in }\Omega\times[0,T]\,,\\
\phi  &= \xi \qquad &&\text{ on }\Gamma\times[0,T]\,,\\
\lim_{x_2\rightarrow-\infty}\nabla \phi&= 0 \qquad &&\text{on }\bT\times[0,T]\,,\\
\xi_t &=-\frac{\varepsilon}{2}\sA^\ell_j\phi,_\ell\sA^k_j\phi,_k- h-\alpha \sA^\ell_2(\sA^k_2\phi,_k),_\ell \nonumber\\
&\quad+\varepsilon \sA^k_2\phi,_k\left(\sA^k_j\phi,_k\tilde{n}_j+\alpha h,_{11}\right)&&\text{ on }\Gamma\times[0,T]\,,\\
h_t&=\sA^k_j\phi,_k\tilde{n}_j+\alpha h,_{11}&&\text{ on }\Gamma\times[0,T]\,,\\
\Delta \delta\psi&= 0  &&\text{ in }\Omega\times[0,T]\,,\\
\delta\psi&= h &&\text{ on }\Gamma\times[0,T]\,,\\
\psi&= e+\delta\psi &&\text{ in }\Omega\times[0,T]\,,\end{align}
\end{subequations}
with the initial data
\begin{align*}
\xi &=\xi_0&&\text{ on }\bT\times\{0\}\,,\\
h&=h_0&&\text{ on }\bT\times\{0\}\,,\\
\end{align*}

\section{Proof of Theorem \ref{theorem1}: Local well-posedness in Sobolev spaces}
\subsection{A smooth approximation of the ALE formulation of the damped water waves}
We consider $\kappa, \delta>0$ two positive constants. The system \eqref{eq:ALE} has a very fine structure that our regularizing procedure has to maintain in order energy estimates can be performed. We define the following smooth approximation of \eqref{eq:all2dimensionless}
\begin{subequations}\label{eq:smoothALE}
\begin{align}
\sA^\ell_j(\sA^k_j\phi,_k),_\ell&=0&&\text{ in }\Omega\times[0,T]\,,\\
\phi  &= \xi^{\kappa} \qquad &&\text{ on }\Gamma\times[0,T]\,,\\
\lim_{x_2\rightarrow-\infty} \nabla \phi&= 0 \qquad &&\text{on }\bT\times[0,T],\\
\xi_t &=\jeps\bigg{(}-\frac{\varepsilon}{2}\sA^\ell_j\phi,_\ell\sA^k_j\phi,_k- h\nonumber\\
&\quad-\alpha \sA^\ell_2(\sA^k_2\phi,_k),_\ell\nonumber\\
&\quad+\varepsilon \sA^k_2\phi,_k\left(\sA^k_j\phi,_k\tilde{n}_j+\alpha h^{\kappa},_{11}\right)\bigg{)}&&\text{ on }\Gamma\times[0,T]\,,\\
h_t&=\jeps\left(\sA^k_j\phi,_k\tilde{n}_j\right)+\alpha h^{\kappa\kappa},_{11}&&\text{ on }\Gamma\times[0,T]\,,\\
\Delta \delta\psi&= 0  &&\text{ in }\Omega\times[0,T]\,,\\
\delta\psi&= h^{\kappa} &&\text{ on }\Gamma\times[0,T]\,,\\
\psi&= e+\delta\psi &&\text{ in }\Omega\times[0,T]\,,
\end{align}
\end{subequations}
where
\begin{align*}
\xi &=\xi^{\delta}_0&&\text{ on }\bT\times\{0\}\,,\\
h&=h^{\delta}_0&&\text{ on }\bT\times\{0\}\,,\\
\end{align*}
are the initial data
and
$$
\tilde{n}=(-\varepsilon h^{\kappa},_1,1).
$$
We observe that we mollify twice in order we have a symmetric regularizing operator. First we need to regularize the trace of the diffeomorphism $\psi$ and also the velocity potential $\phi$. After this regularizing procedure, we are able to solve the elliptic problems associated to the velocity and the diffeomorphism. Furthermore, due to these mollifiers, the high order term in the equation for $\xi$ is
$$
\sA^2_2(\sA^2_2\phi,_2),_2\approx \xi^{\kappa},_{11}+\text{ nonlinear terms}.
$$
Then the parabolic gain of regularity translates into a negative contribution of the type
$$
-\int_\bT (\xi^\kappa,_1)^2dx.
$$
Then, in order the structure of the system remains unchanged and the energy estimates can be performed, we have to mollify also the high order terms $h^{\kappa\kappa},_{11}$ present in the equation for $h$ and $\xi$. This is due to the fact that the system \eqref{eq:ALE} is of cross-diffusion type and, without considering smoothing operators, we will have to bound terms akin to 
$$
\int_{\bT} \xi\sA^k_2\phi,_k h,_{11}dx.
$$
Thus, we need to have $h^{\kappa\kappa},_{11}$ in order we can use the previous negative contribution to avoid loss of derivatives.

\subsection{Existence for the regularized damped water waves system}
To prove the existence of local solution for the system \eqref{eq:smoothALE}, we use a standard fixed point scheme (see \cite{CGS}). We define the Banach space
$$
X=C([0,T],H^{3}(\bT))\times C([0,T],H^{3}(\bT))
$$
and consider a pair
$$
(\bar{h},\bar{\xi})\in X.
$$
Then we define the following \emph{linear} problem:
\begin{subequations}\label{eq:smoothALElin}
\begin{align}
\bar{\sA}^\ell_j(\bar{\sA}^k_j\phi,_k),_\ell&=0&&\text{ in }\Omega\times[0,T]\,,\\
\phi  &= \bar{\xi}^{\kappa} \qquad &&\text{ on }\Gamma\times[0,T]\,,\\
\lim_{x_2\rightarrow-\infty}\nabla \phi &= 0 \qquad &&\text{on }\bT\times[0,T],\\
\xi_t &=\jeps\bigg{(}-\frac{\varepsilon}{2}\bar{\sA}^\ell_j\phi,_\ell \bar{\sA}^k_j\phi,_k- h\nonumber\\
&\quad-\alpha \bar{\sA}^\ell_2(\bar{\sA}^k_2\phi,_k),_\ell\nonumber\\
&\quad+\varepsilon \bar{\sA}^k_2\phi,_k\left(\bar{\sA}^k_j\phi,_k\tilde{n}_j+\alpha h^{\kappa},_{11}\right)\bigg{)}&&\text{ on }\Gamma\times[0,T]\,,\\
h_t&=\jeps\left(\bar{\sA}^k_j\phi,_k\tilde{n}_j\right)+\alpha h^{\kappa\kappa},_{11}&&\text{ on }\Gamma\times[0,T]\,,\\
\Delta \delta\psi&= 0  &&\text{ in }\Omega\times[0,T]\,,\\
\delta\psi&= h^{\kappa} &&\text{ on }\Gamma\times[0,T]\,,\\
\psi&= e+\delta\psi &&\text{ in }\Omega\times[0,T]\,,
\end{align}
\end{subequations}
with initial data
\begin{align*}
\xi &=\xi^{\delta}_0&&\text{ on }\bT\times\{0\}\,,\\
h&=h^{\delta}_0&&\text{ on }\bT\times\{0\}\,,\\
\end{align*}

Now we consider the mapping
$$
S[\bar{h},\bar{\xi}]=(h,\xi)
$$
It is easy to see that 
$$
S:X\rightarrow X.
$$
and, for a short enough time interval $T=T(\kappa)$, the mapping is contractive. Thus, there exists a fixed point that is our approximate solution for \eqref{eq:smoothALE}. 
\subsection{The energy functional}
We observe that the solution to the approximate problem satisfies
$$
(h,\xi)\in C([0,T],H^{3})\times C([0,T],H^3).
$$
Furthermore, 
$$
h^\kappa \in L^2(0,T;H^{4})\;,\;\xi^\kappa\in L^2(0,T;H^{4}).
$$
The purpose of this section is to provide with $\kappa-$independent estimates that ensure a lower bound 
$$
0<T_*\leq T(\kappa).
$$
In particular, we want to bound the following energy
$$
\mathscr{E}(t)=\max_{0\leq s\leq t}\left(|h(s)|_{3}^2+|\xi(s)|_{3}^2\right)+\int_0^t\|\nabla\phi(s)\|_{2.5}^2ds.
$$

By restricting $T(\kappa)$ if necessary, we can ensure that
\begin{enumerate}
\item 
$$
\mathscr{E}(t)< z^*, \text{ for all }0\leq t\leq T(\kappa),
$$
for 
$$
z^*=2\mathscr{E}(0)=2(|h_0|_{3}^2+|\xi_0|_{3}^2),
$$ 
\item 
$$
\max\{|h(t)|_2,|\xi(t)|_2\}< \gamma<1, \text{ for all }0\leq t\leq T(\kappa),
$$
\item for a fixed $0<\beta\ll1$,
\begin{equation}\label{bootstrap2}
\|\sA(t)-\text{Id}\|_{L^\infty}< \beta, \text{ for all }0\leq t\leq T(\kappa),
\end{equation}
where $\text{Id}=\delta^i_k$ denotes the identity matrix.
\end{enumerate}

We want to find a polynomial estimate of the form
$$
\mathscr{E}(t)\leq \mathcal{M}_0+t^\sigma\mathcal{Q}(\mathscr{E}(t)),
$$
for certain polynomial $\mathcal{Q}$ and positive constants $\mathcal{M}_0$ and $\sigma$. Such an inequality implies the existence of the desired $T_*$ and the uniform-in-$\kappa$ bound (see \cite{coutand2006interaction} for further details)
$$
\mathscr{E}(t)\leq 2\mathcal{M}_0\text{ for all }0\leq t\leq T_*.
$$
We take $T(\kappa)$ small enough so
\begin{equation}\label{smalltime}
\mathcal{Q}(10\mathscr{E}(0))T(\kappa)^{1/8}\leq \lambda,
\end{equation}
for certain $\lambda<\beta\ll1$  fixed such that
$$
0<\lambda<\gamma-\max\{|h_0|_2,|\xi_0|_2\}.
$$

\subsection{$\kappa-$independent estimates}
\subsubsection*{Estimates for the diffeomorphism $\psi$}
Using standard elliptic estimates, we find that
$$
\|\nabla\delta\psi(t)\|_{2.5}\leq C|h(t)|_{3}.
$$
We also have
$$
\|\nabla\psi(t)-\text{Id}\|_{2.5}\leq C|h(t)|_{3}.
$$
In general,
\begin{equation}\label{estimatepsi}
\|\nabla\psi-\text{Id}\|_{r}\leq C|h^{\kappa}|_{r+0.5},\forall\,0<r\leq2.5.
\end{equation}

\subsubsection*{Smallness in lower norms for the interface}
Using the equation for $h_t$, we have the bound
$$
|h_t|_{L^2}\leq C|\sA|_{L^\infty}|\nabla\phi|_{L^\infty}(1+|h|_1)+|h|_2,
$$
so, using Sobolev embedding
$$
H^{1}(\bT)\subset L^\infty(\bT),
$$ 
together with
$$
|J^{-1}|=\left|\frac{1}{1+\delta\psi,_2}\right|\leq \frac{1}{1-C|h^{\kappa}|_2},
$$
the estimate
$$
|\sA|_{L^\infty}=|J^{-1}J\sA|_{L^\infty}\leq \frac{|J\sA|_{L^\infty}}{1-C|h^{\kappa}|_2}\leq C|\nabla\psi-\text{Id}+\text{Id}|_{L^\infty}\leq C(1+|\nabla\psi-\text{Id}|_{L^\infty}),
$$
and the Trace theorem for $H^s(\Gamma)$, we find that
\begin{align*}
\int_0^t|h_t(s)|_{L^2}^2ds&\leq C\int_0^t (1+\|\nabla\psi(s)-\text{Id}\|_{1.5}^2)\|\nabla\phi(s)\|_{1.5}^2(1+|h(s)|_1)^2+|h(s)|_2^2ds\\
&\leq C\mathscr{E}(t)(1+\mathscr{E}(t))^2+tC\mathscr{E}(t).
\end{align*}
As a consequence, we obtain that
\begin{align*}
|h(t)-h^\delta_0|_{L^2}&\leq \sqrt{t}\left(\int_0^t|h_t(s)|_{L^2}^2ds\right)^{1/2}\\
&\leq \sqrt{t}\left(C\mathscr{E}(t)(1+\mathscr{E}(t))^2+tC\mathscr{E}(t)\right)^{0.5}\\
&\leq \sqrt{t}\left(Cz^*(1+z^*)^2+tCz^*\right)^{1/2} .
\end{align*}
Taking a sufficiently small time $T(\kappa)$ as in \eqref{smalltime}, we find that
$$
|h(t)-h^\delta_0|_{2}=|h(t)-h^\delta_0|^{1/3}_{L^2}|h(t)-h^\delta_0|_{3}^{2/3}\leq t^{1/6}\left(z^*(1+z^*)^2\right)^{1/6}C(z^*)^{2/3} \leq \lambda,
$$
for a small enough $\lambda$. Taking
$$
0<\lambda<\gamma-\max\{|h_0|_2,|\xi_0|_2\}.
$$
small enough, this smallness property has the following important consequence
\begin{equation}\label{normconstant}
|h(t)|_{2}\leq |h^\delta_0|_{2}+\lambda<\gamma,
\end{equation}
in the $\kappa-$dependent time interval $[0,T(\kappa)]$.

Then, we also estimate
$$
\|J(t)-1\|_{1.5}\leq \|\nabla\delta\psi\|_{1.5}\leq |h^{\kappa}|_2, 
$$
which implies that
$$
0<C\leq J(t)\leq C^{-1}
$$
in the $\kappa-$dependent time interval $[0,T(\kappa)]$.

\subsubsection*{Smallness in lower norms for the trace of the potential}
Using now the equation for $\xi_t$ together with the fact that $\sA^{1}_2=0$, we find that
\begin{align*}
|\xi_t|_{L^2}&\leq C\bigg{(}|\sA|_{L^\infty}^2|\nabla\phi|_{L^\infty}^2 +|h|_{L^2}+ |\sA|_{L^\infty}|\sA,_2|_{L^2}|\nabla\phi|_{L^\infty}+|\sA|_{L^\infty}^2|\nabla \phi,_2|_{L^2}\\
&\quad+ |\sA|_{L^\infty}|\nabla\phi|_{L^\infty}\left(|\sA|_{L^\infty}|\nabla\phi|_{L^\infty}(1+|h|_1)+|h|_2\right)\bigg{)}.
\end{align*}
Since
$$
|\sA,_2|_{L^2}\leq C|\sA,_2|_{0.25}\leq C\|\sA,_2\|_{0.75}\leq C\|\sA,_2\|_{1}, 
$$
and invoking the Trace Theorem together with Sobolev embedding, we have that
\begin{align}
|\xi_t|_{L^2}&\leq C\bigg{(}(1+\|\nabla\psi-\text{Id}\|_{1.5})^2\|\nabla\phi\|_{1.5}^2 +|h|_{L^2}+ (1+\|\nabla\psi-\text{Id}\|_{1.5})\|\nabla\delta\psi\|_{2}\|\nabla\phi\|_{1.5}\nonumber\\
&\quad+(1+\|\nabla\psi-\text{Id}\|_{1.5})^2\|\nabla \phi\|_{2}\label{estimatexit}\\
&\quad+ (1+\|\nabla\psi-\text{Id}\|_{1.5})\|\nabla\phi\|_{1.5}\left((1+\|\nabla\psi-\text{Id}\|_{1.5})\|\nabla\phi\|_{1.5}(1+|h|_1)+|h|_2\right)\bigg{)}\nonumber.
\end{align}
As a consequence
$$
\int_0^t|\xi_t(s)|_{L^2}^2ds\leq t\mathcal{P}(\mathscr{E}(t)),
$$
for certain (explicit) polynomial $\mathcal{P}$. Thus, we find that
$$
|\xi(t)-\xi^\delta_0|_{L^2}\leq \sqrt{t}\left(\int_0^t|\xi_t(s)|_{L^2}^2ds\right)^{1/2}\leq \sqrt{t}\mathcal{P}(\mathscr{E}(t)).
$$
Taking a sufficiently small time $T(\kappa)$ as in \eqref{smalltime}, we find that
$$
|\xi(t)-\xi^\delta_0|_{2}=|\xi(t)-\xi^\delta_0|^{1/3}_{L^2}|\xi(t)-\xi^\delta_0|_{3}^{2/3}\leq \lambda,
$$
for a small enough $\lambda$. Taking $\lambda=\lambda(h_0,\gamma)$ small enough as before, we conclude
\begin{equation}\label{normconstant2}
|\xi(t)|_{2}\leq |\xi^\delta_0|_{2}+\lambda<\gamma,\;\forall\,t\in [0,T(\kappa)].
\end{equation}

\subsubsection*{Estimates for the velocity potential}
Using Piola's identity, we have that the velocity potential satisfies
\begin{subequations}\label{eq:potential0}
\begin{align}
(J\sA^\ell_j\sA^k_j\phi,_k),_\ell&=0&&\text{ in }\Omega\times[0,T]\,,\\
\phi  &= \xi^{\kappa} \qquad &&\text{ on }\Gamma\times[0,T]\,,\\
\lim_{x_2\rightarrow-\infty}\nabla \phi &= 0 \qquad &&\text{on }\bT\times[0,T]\,.
\end{align}
\end{subequations}
Equivalently, we can write the previous system as
\begin{subequations}\label{eq:potential}
\begin{align}
\Delta\phi&=-\nabla\cdot\left(J\sA\sA^\intercal-\text{Id})\nabla\phi\right)&&\text{ in }\Omega\times[0,T]\,,\\
\phi  &= \xi^{\kappa} \qquad &&\text{ on }\Gamma\times[0,T]\,,\\
\lim_{x_2\rightarrow-\infty}\nabla \phi &= 0 \qquad &&\text{on }\bT\times[0,T]\,.
\end{align}
\end{subequations}
Let us decompose $ \phi $ as
\begin{equation}\label{phidecomposition}
\phi = \phi_1 +\phi_2 , 
\end{equation}
where
\begin{align}\label{eq:eq_phi1}
& \left\lbrace
\begin{aligned}
& \Delta \phi_1 =0, & \text{ in }\Omega\times[0,T], \\
& \phi_1 = \xi^{\kappa}, & \text{ on } \Gamma\times[0,T], 
\end{aligned}
\right. 
\end{align}
and
\begin{align}\label{eq:eq_phi2}
\left\lbrace
\begin{aligned}
& \Delta\phi_2 = - \nabla \cdot \bra{\pare{ J \sA \sA^\intercal -\text{Id}  } \nabla( \phi_1+\phi_2)\big. }, & \text{ in }\Omega\times[0,T], \\
& \phi_2 = 0, & \text{ on } \Gamma\times[0,T].
\end{aligned}
\right. 
\end{align}
We compute
\begin{equation}\label{matrixelliptic}
J \sA \sA^\intercal -\text{Id}=\left(\begin{array}{cc} \delta\psi,_2& -\delta\psi,_1\\ -\delta\psi,_1& \frac{\displaystyle(\delta\psi,_1)^2}{\displaystyle1+\delta\psi,_2}-\frac{\displaystyle\delta\psi,_2}{\displaystyle 1+\delta\psi,_2}\end{array}\right)
\end{equation}
We observe that
\begin{equation}\label{eq:phi1}
\phi_1\pare{x_1, x_2, t} = e^{x_2\Lambda}\xi^{\kappa}\pare{x_1, t},
\end{equation}
contains the linear contribution while $\phi_2$ is purely nonlinear. We have that
\begin{equation}\label{estimatephi1}
\|\nabla\phi_1\|_{r}\leq C|\xi^{\kappa}|_{r+0.5},\forall\,0<r\leq2.5.
\end{equation}
Let us first consider the estimate for a lower order norm of $\phi_2$. Elliptic estimates (see \cite{Lannes2013,CGS,granero2019well} for further details on similar problems), the Banach Algebra property of the Sobolev spaces with enough regularity, \eqref{estimatepsi} and \eqref{matrixelliptic} then give us
\begin{align*}
\|\nabla \phi_2\|_{1.5}&\leq \|\pare{ J \sA \sA^\intercal -\text{Id}  } \nabla( \phi_1+\phi_2)\|_{1.5}\\
&\leq \|\pare{ J \sA \sA^\intercal -\text{Id}  } \nabla\phi_1\|_{1.5}+\|\pare{ J \sA \sA^\intercal -\text{Id}  } \nabla\phi_2\|_{1.5}\\
&\leq C\| J \sA \sA^\intercal -\text{Id}  \|_{1.5} \|\nabla\phi_1\|_{1.5}+C\| J \sA \sA^\intercal -\text{Id}  \|_{1.5}\|\nabla\phi_2\|_{1.5}\\
&\leq C\|\nabla\delta\psi\|_{1.5}|\xi^{\kappa}|_{2}+C\|\nabla\delta\psi\|_{1.5}\|\nabla\phi_2\|_{1.5}\\
&\leq C|h^{\kappa}|_2|\xi^{\kappa}|_{2}+C|h^{\kappa}|_2\|\nabla\phi_2\|_{1.5}.
\end{align*}
Using the smallness of $|h^{\kappa}|_2$, we can absorb the last term into the left hand side and we find that
\begin{align}\label{estimatephi2LOWER}
\|\nabla \phi_2\|_{1.5}&\leq C|\xi^{\kappa}|_{2}.
\end{align}
The higher order norm can be bounded similarly,
\begin{align*}
\|\nabla \phi_2\|_{2.5}&\leq \|\pare{ J \sA \sA^\intercal -\text{Id}  } \nabla( \phi_1+\phi_2)\|_{2.5}\\
&\leq \|\pare{ J \sA \sA^\intercal -\text{Id}  } \nabla\phi_1\|_{2.5}+\|\pare{ J \sA \sA^\intercal -\text{Id}  } \nabla\phi_2\|_{2.5}\\
&\leq C\left(\| J \sA \sA^\intercal -\text{Id}  \|_{2.5} \|\nabla\phi_1\|_{L^\infty}+\| J \sA \sA^\intercal -\text{Id}  \|_{L^\infty} \|\nabla\phi_1\|_{2.5}\right.\\
&\quad+\left.\|J \sA \sA^\intercal -\text{Id}  \|_{L^\infty} \|\nabla\phi_2\|_{2.5}+\| J \sA \sA^\intercal -\text{Id} \|_{2.5} \|\nabla\phi_2\|_{L^\infty}\right).
\end{align*}
Then, using the Sobolev embedding together with the previous bounds for $\|\nabla\phi_2\|_{1.5}$, $\|\nabla\phi_1\|_{r}$ and $\|\nabla\delta\psi\|_{r}$, we find that
\begin{align*}
\|\nabla \phi_2\|_{2.5}&\leq C\left(\| J \sA \sA^\intercal -\text{Id}  \|_{2.5} \|\nabla\phi_1\|_{1.5}+\| J \sA \sA^\intercal -\text{Id}  \|_{1.5} \|\nabla\phi_1\|_{2.5}\right.\\
&\quad+\left.\|J \sA \sA^\intercal -\text{Id}  \|_{1.5} \|\nabla\phi_2\|_{2.5}+\| J \sA \sA^\intercal -\text{Id} \|_{2.5} \|\nabla\phi_2\|_{1.5}\right)\\
&\leq C\left(|h^{\kappa} |_{3} |\xi^{\kappa}|_{2}+|h^{\kappa} |_{2} |\xi^{\kappa}|_{3}+|h^{\kappa}|_{2} \|\nabla\phi_2\|_{2.5}\right).
\end{align*}
Thus, using again the smallness of $|\xi^{\kappa}|_2$ and $|h^{\kappa}|_2$, we conclude that
\begin{equation}\label{estimatephi2}
\|\nabla \phi_2\|_{2.5}\leq C\left(|h^{\kappa} |_{3}+|\xi^{\kappa}|_{3}\right).
\end{equation}
Finally, we also have that
\begin{align*}
\|\nabla \phi_2\|_{3.5}&\leq \|\pare{ J \sA \sA^\intercal -\text{Id}  } \nabla( \phi_1+\phi_2)\|_{3.5}\\
&\leq \|\pare{ J \sA \sA^\intercal -\text{Id}  } \nabla\phi_1\|_{3.5}+\|\pare{ J \sA \sA^\intercal -\text{Id}  } \nabla\phi_2\|_{3.5}\\
&\leq C\left(\| J \sA \sA^\intercal -\text{Id}  \|_{3.5} \|\nabla\phi_1\|_{L^\infty}+\| J \sA \sA^\intercal -\text{Id}  \|_{L^\infty} \|\nabla\phi_1\|_{3.5}\right.\\
&\quad+\left.\|J \sA \sA^\intercal -\text{Id}  \|_{L^\infty} \|\nabla\phi_2\|_{3.5}+\| J \sA \sA^\intercal -\text{Id} \|_{3.5} \|\nabla\phi_2\|_{L^\infty}\right)\\
&\leq C\left(\| J \sA \sA^\intercal -\text{Id}  \|_{3.5} \|\nabla\phi_1\|_{1.5}+\| J \sA \sA^\intercal -\text{Id}  \|_{1.5} \|\nabla\phi_1\|_{3.5}\right.\\
&\quad+\left.\|J \sA \sA^\intercal -\text{Id}  \|_{1.5} \|\nabla\phi_2\|_{3.5}+\| J \sA \sA^\intercal -\text{Id} \|_{3.5} \|\nabla\phi_2\|_{1.5}\right)\\
&\leq C\left(\| \nabla\delta\psi  \|_{3.5}|\xi^\kappa|_{2}+\| \nabla\delta\psi \|_{1.5} |\xi^\kappa|_{4}\right.\\
&\quad+\left.\| \nabla\delta\psi \|_{1.5} \|\nabla\phi_2\|_{3.5}+\| \nabla\delta\psi \|_{3.5} |\xi^\kappa|_{2}\right)\\
&\leq C\left(|h^\kappa|_{4}|\xi^\kappa|_{2}+| h^\kappa |_{2} |\xi^\kappa|_{4}\right.\\
&\quad+\left.| h^\kappa|_{2} \|\nabla\phi_2\|_{3.5}+| h^\kappa |_{4} |\xi^\kappa|_{2}\right).
\end{align*}
Using the smallness of $|h|_2$, we conclude that
\begin{equation}\label{estimatephi2B}
\|\nabla \phi_2\|_{3.5}\leq C\left(|h^\kappa|_{4}|\xi^\kappa|_{2}+| h^\kappa |_{2} |\xi^\kappa|_{4}+| h^\kappa |_{4} |\xi^\kappa|_{2}\right).
\end{equation}

\subsubsection*{Estimates for the interface}
We observe that the mean is preserved by the equation for $h$. In other words,
$$
\int_{\bT}h(x_1,t)dx_1=0\;\forall\,0\leq t.
$$
As a consequence, 
$$
|h|_3\leq C|h,_{111}|_{L^2}.
$$
Testing the equation for the interface against $\Lambda^{6}h$, using the symmetry of the smoothing operator and integrating by parts, we find that
\begin{align*}
\frac{1}{2}\frac{d}{dt}|h,_{111}|_{L^2}^2&=\int_\bT\sA^k_j\phi,_k\tilde{n}_j\Lambda^6h^\kappa dx_1 -|h^\kappa,_{1111}|^2_{L^2}\\
&=-\int_\bT\sA^k_j\phi,_k\tilde{n}_j\partial_1^6h^\kappa dx_1 -|h^\kappa,_{1111}|^2_{L^2}\\
&=-\int_\bT\partial_1^2(\sA^k_j\phi,_k\tilde{n}_j)\partial_1^4h^\kappa dx_1 -|h^\kappa,_{1111}|^2_{L^2}.
\end{align*}
Using
$$
\psi,_2(x_1,0,t)=\Lambda h^{\kappa},
$$
we compute
$$
\sA^k_j\phi,_k\tilde{n}_j=-\sA^k_1\phi,_k h^{\kappa},_1+\sA^k_2\phi,_k=-\xi^{\kappa},_1 h^{\kappa},_1+\frac{(h^{\kappa},_1)^2}{1+\Lambda h^{\kappa}}\phi,_2 +\frac{\phi,_2}{1+\Lambda h^{\kappa}}.
$$
We note that, using the smallness of $h$ in $H^2$ \eqref{normconstant},
$$
\left|\frac{1}{1+\Lambda h^{\kappa}}\right|_{L^\infty}\leq \frac{1}{1-|\Lambda h^{\kappa}|_{L^\infty}}\leq \frac{1}{1-|h|_{2}}\leq 2.
$$
Using the Sobolev embedding and interpolation, we have that
$$
|\sA^k_j\phi,_k\tilde{n}_j|_{2}\leq C\left(|\xi|_3|h|_{2}+|\xi|_{2}|h|_{3}+|h|_2^2|\phi,_2|_{2}+|h|_2|\phi,_2|_{1}|h|_3+|\phi,_2|_{2}+|\phi,_2|_{1}|h|_3\right).
$$
Invoking the elliptic estimates of the previous section \eqref{estimatephi1} and \eqref{estimatephi2}, together with the Trace theorem, we find that
$$
|\phi,_2|_{2}\leq C\|\nabla \phi\|_{2.5}\leq C(|h|_3+|\xi|_3).
$$
As a consequence,
$$
|\sA^k_j\phi,_k\tilde{n}_j|_{2}\leq C\sqrt{\mathscr{E}(t)}.
$$
Then, using Young's inequality, we conclude that
\begin{align}\label{estimateh}
\frac{d}{dt}|h,_{111}|_{L^2}^2&\leq C\mathscr{E}(t) -|h^\kappa,_{1111}|^2_{L^2}.
\end{align}
\subsubsection*{Estimates for the trace of the velocity potential}
Since \eqref{estimatexit}, the mean satisfies the following ordinary differential inequality
\begin{align*}\frac{d}{dt}\langle \xi(t)\rangle=\frac{d}{dt}\int_\bT\xi(x_1,t) dx_1 &= \int_{\bT}\xi_t(x_1,t)dx_1\leq |\xi_t|_{L^1}\leq \mathcal{P}(\mathscr{E}(t))\,.
\end{align*}
Now we have that 
$$
|\xi|_{3}\leq |\xi-\langle \xi\rangle+\langle \xi\rangle|_{3}\leq |\xi-\langle \xi\rangle|_3+|\langle \xi\rangle|_{L^2}\leq |\xi-\langle \xi\rangle|_3+C|\langle \xi\rangle|\leq |\xi,_{111}|_{L^2}+t\mathcal{P}(\mathscr{E}(t)).
$$
Thus, it is enough to estimate $|\xi,_{111}|_{L^2}$.

Let us now consider the term $ \sA_2^2 \partial_2 \pare{\sA_2^2 \partial_2 \phi} $, which is the higher-order nonlinear term in the equation for $\xi$. We want to extract the linear contribution. This linear contribution has a sign and the consequent dissipative effect will play an important role in the estimates. We have that
\begin{align*}
\sA_2^2 \pare{\sA_2^2 \phi,_2},_2 
% = \pare{\sA_2^2}^2 \partial_2^2 \phi + \sA_2^2 \partial_2 \sA_2^2 \partial_2 \phi, \\
%%----------------------------------------------------------
& = \frac{1}{\pare{1+\delta\psi,_2}^2}  \ \phi,_{22} - \frac{\delta\psi,_{22}}{\pare{1+\delta\psi,_2}^3} \phi,_2, \\
%%----------------------------------------------------------
& = \frac{1}{\pare{1+\delta\psi,_2}^2}  \ \pare{ \phi_1+\phi_2 },_{22} - \frac{\delta\psi,_{22}}{\pare{1+\delta\psi,_2}^3} \pare{ \phi_1+\phi_2 },_2. 
\end{align*}
Using the explicit expression of $\phi_1$ and $\psi$, we compute
$$
\delta\psi,_{22}(x_1,0,t)=-h^{\kappa},_{11},
$$
$$
\phi_1,_2(x_1,0,t)=\Lambda \xi^{\kappa},
$$
$$
\phi_1,_{22}(x_1,0,t)=-\xi^{\kappa},_{11}
$$
\begin{align*}
\left[ \frac{1}{\pare{1+\delta\psi,_2}^2}  \ \phi_1,_{22} - \frac{\delta\psi,_{22}}{\pare{1+\delta\psi,_2}^3}\phi_1,_2\right]\bigg{|}_{x_2=0} & = -\frac{\xi^{\kappa},_{11} }{\pare{1+\Lambda h^{\kappa}}^2}  + \frac{h^{\kappa},_{11}}{\pare{1+\Lambda h^{\kappa}}^3} \Lambda \xi^{\kappa}, \\
& = -\xi^{\kappa},_{11} + \frac{2\Lambda h^{\kappa} + \pare{\Lambda h^{\kappa} }^2}{\pare{1+\Lambda h^{\kappa}}^2} \ \xi^{\kappa},_{11} + \frac{h^{\kappa},_{11}}{\pare{1+\Lambda h^{\kappa}}^3} \Lambda \xi^{\kappa}.
\end{align*}
As a consequence
\begin{align}
-\sA^\ell_2(\sA^k_2\phi,_k),_\ell\bigg{|}_{x_2=0}&=-\sA^2_2(\sA^2_2\phi,_2),_2\bigg{|}_{x_2=0}\nonumber\\
&=\frac{-1}{\pare{1+\Lambda h^\kappa}^2}  \ \phi_2,_{22}\bigg{|}_{x_2=0} - \frac{h^\kappa,_{11}}{\pare{1+\Lambda h^\kappa}^3} \phi_2,_2\bigg{|}_{x_2=0}\nonumber\\
&\quad +\xi^{\kappa},_{11} - \frac{2\Lambda h^{\kappa} + \pare{\Lambda h^{\kappa} }^2}{\pare{1+\Lambda h^{\kappa}}^2} \ \xi^{\kappa},_{11} - \frac{h^{\kappa},_{11}}{\pare{1+\Lambda h^{\kappa}}^3} \Lambda \xi^{\kappa}.\label{eq:aux1}
\end{align}
Thus, the equation for $\xi_t$ reads
\begin{align*}
\xi_t &=\jeps\bigg{(}-\frac{1}{2}\sA^\ell_j\phi,_\ell\sA^k_j\phi,_k- h-\frac{\phi_2,_{22}}{\pare{1+\Lambda h^\kappa}^2}   - \frac{h^\kappa,_{11}\phi_2,_2}{\pare{1+\Lambda h^\kappa}^3}+\xi^{\kappa},_{11} \nonumber\\
&\quad  - \frac{2\Lambda h^{\kappa} + \pare{\Lambda h^{\kappa} }^2}{\pare{1+\Lambda h^{\kappa}}^2} \ \xi^{\kappa},_{11} - \frac{h^{\kappa},_{11}\Lambda \xi^{\kappa}}{\pare{1+\Lambda h^{\kappa}}^3} + \frac{\phi,_2}{1+\Lambda h^\kappa}\left(\sA^k_j\phi,_k\tilde{n}_j+ h^{\kappa},_{11}\right)\bigg{)}&&\text{ on }\Gamma\times[0,T]\,.
\end{align*}
Then we have that
\begin{align*}
\frac{d}{dt}|\xi,_{111}|_{L^2}^2&=-|\xi^\kappa,_{1111}|_{L^2}^2-\int_{\bT}\xi^\kappa,_{1111}\bigg{(}-\frac{1}{2}\sA^\ell_j\phi,_\ell\sA^k_j\phi,_k- h-\frac{\phi_2,_{22}}{\pare{1+\Lambda h^\kappa}^2}   - \frac{h^\kappa,_{11}\phi_2,_2}{\pare{1+\Lambda h^\kappa}^3} \nonumber\\
&\quad  - \frac{2\Lambda h^{\kappa} + \pare{\Lambda h^{\kappa} }^2}{\pare{1+\Lambda h^{\kappa}}^2} \ \xi^{\kappa},_{11} - \frac{h^{\kappa},_{11}\Lambda \xi^{\kappa}}{\pare{1+\Lambda h^{\kappa}}^3} + \frac{\phi,_2}{1+\Lambda h^\kappa}\left(\sA^k_j\phi,_k\tilde{n}_j+ h^{\kappa},_{11}\right)\bigg{)},_{11}dx\\
&=-|\xi^\kappa,_{1111}|_{L^2}^2+ I_1+ I_2+ I_3+ I_4+ I_5+ I_6+ I_7+I_8.
\end{align*}
with
\begin{align*}
I_1&=\int_{\bT}\xi^\kappa,_{1111}\bigg{(}\frac{1}{2}\sA^\ell_j\phi,_\ell\sA^k_j\phi,_k\bigg{)},_{11}dx,\\
I_2&=\int_{\bT}\xi^\kappa,_{1111}h,_{11}dx,\\
I_3&=\int_{\bT}\xi^\kappa,_{1111}\bigg{(}\frac{\phi_2,_{22}}{\pare{1+\Lambda h^\kappa}^2}\bigg{)},_{11}dx,\\
I_4&=\int_{\bT}\xi^\kappa,_{1111}\bigg{(}\frac{h^\kappa,_{11}\phi_2,_2}{\pare{1+\Lambda h^\kappa}^3}\bigg{)},_{11}dx,\\
I_5&=\int_{\bT}\xi^\kappa,_{1111}\bigg{(}\frac{2\Lambda h^{\kappa} + \pare{\Lambda h^{\kappa} }^2}{\pare{1+\Lambda h^{\kappa}}^2} \xi^{\kappa},_{11}\bigg{)},_{11}dx,\\
I_6&=\int_{\bT}\xi^\kappa,_{1111}\bigg{(}\frac{h^{\kappa},_{11}\Lambda \xi^{\kappa}}{\pare{1+\Lambda h^{\kappa}}^3}\bigg{)},_{11}dx,\\
I_7&=-\int_{\bT}\xi^\kappa,_{1111} \bigg{(}\frac{\phi,_2}{1+\Lambda h^\kappa}\left(\sA^k_j\phi,_k\tilde{n}_j\right)\bigg{)},_{11}dx,\\
I_8&=-\int_{\bT}\xi^\kappa,_{1111} \bigg{(}\frac{\phi,_2}{1+\Lambda h^\kappa} h^{\kappa},_{11}\bigg{)},_{11}dx.
\end{align*}
Using the elliptic estimates \eqref{estimatephi1} and \eqref{estimatephi2}, we find the bound
$$
I_1\leq |\xi^\kappa,_{1111}|_{L^2}C(|h|_3^2+|\xi|_3^2).
$$
Similarly, using again the elliptic estimates \eqref{estimatephi1} and \eqref{estimatephi2} together with the smallness of $|h|_2+|\xi|_2$, we obtain that
$$
I_7\leq C|\xi^\kappa,_{1111}|_{L^2}(|h|_3^2+|\xi|_3^2).
$$
Applying Cauchy-Schwarz inequality we obtain that
$$
I_ 2\leq |\xi^\kappa,_{1111}|_{L^2}|h|_2.
$$
To bound the terms $I_4$, $I_5$, $I_6$ and $I_8$ we need to use the smallness of $|h|_2+|\xi|_2$. This is due to the cross-diffusion character of the system. It is easy to observe that
\begin{align*}
I_4&\leq\int_{\bT}\xi^\kappa,_{1111}\frac{h^\kappa,_{1111}\phi_2,_2}{\pare{1+\Lambda h^\kappa}^3}dx+C|\xi^\kappa,_{1111}|_{L^2}(|h|_3^2+|\xi|_3^2)\\
&\leq C|\xi^\kappa,_{1111}|_{L^2}|h^\kappa,_{1111}|_{L^2}\|\nabla \phi_2\|_{L^\infty}+C|\xi^\kappa,_{1111}|_{L^2}(|h|_3^2+|\xi|_3^2)\\
&\leq C|\xi^\kappa,_{1111}|_{L^2}|h^\kappa,_{1111}|_{L^2}\|\nabla \phi_2\|_{1.5}+C|\xi^\kappa,_{1111}|_{L^2}(|h|_3^2+|\xi|_3^2)\\
&\leq C|\xi^\kappa,_{1111}|_{L^2}|h^\kappa,_{1111}|_{L^2}|\xi^{\kappa}|_{2}+C|\xi^\kappa,_{1111}|_{L^2}(|h|_3^2+|\xi|_3^2),
\end{align*}
where we have used \eqref{estimatephi2LOWER}. We can proceed analogously for the terms $I_5, I_6$ and $I_8$. Then we obtain that
\begin{align*}
I_5&\leq C|\xi^\kappa,_{1111}|_{L^2}^2|h^{\kappa}|_{2}+C|\xi^\kappa,_{1111}|_{L^2}(|h|_3^2+|\xi|_3^2),\\
I_6&\leq C|\xi^\kappa,_{1111}|_{L^2}|h^\kappa,_{1111}|_{L^2}|\xi^{\kappa}|_{2}+C|\xi^\kappa,_{1111}|_{L^2}(|h|_3^2+|\xi|_3^2),\\
I_8&\leq C|\xi^\kappa,_{1111}|_{L^2}|h^\kappa,_{1111}|_{L^2}\|\nabla \phi\|_{L^\infty}+C|\xi^\kappa,_{1111}|_{L^2}(|h|_3^2+|\xi|_3^2)\\
&\leq C|\xi^\kappa,_{1111}|_{L^2}|h^\kappa,_{1111}|_{L^2}\|\nabla \phi\|_{1.5}+C|\xi^\kappa,_{1111}|_{L^2}(|h|_3^2+|\xi|_3^2)\\
&\leq C|\xi^\kappa,_{1111}|_{L^2}|h^\kappa,_{1111}|_{L^2}|\xi^{\kappa}|_{2}+C|\xi^\kappa,_{1111}|_{L^2}(|h|_3^2+|\xi|_3^2).
\end{align*}
We are left with $I_3$. The term $I_3$ is challenging because there is a term with four derivatives of $\phi_2$. We have that
$$
I_3=J_1+J_2,
$$
with
$$
J_1=\int_{\bT}\xi^\kappa,_{1111}\frac{\phi_2,_{1122}}{\pare{1+\Lambda h^\kappa}^2}dx,
$$
$$
J_2\leq C|\xi^\kappa,_{1111}|_{L^2}(|h|_3^2+|\xi|_3^2).
$$
Using \eqref{estimatephi2B} together with the Trace Theorem, we find that
\begin{align*}
J_1&\leq |\xi^\kappa,_{1111}|_{L^2}\|\nabla\phi_2,_{2}\|_{2.5}\\
&\leq |\xi^\kappa,_{1111}|_{L^2}\|\nabla\phi_2\|_{3.5}\\
&\leq C|\xi^\kappa,_{1111}|_{L^2}\left(|h^\kappa,_{1111}|_{L^2}|\xi^\kappa|_{2}+| h^\kappa |_{2} |\xi^\kappa,_{1111}|_{L^2}+| h^\kappa,_{1111} |_{L^2} |\xi^\kappa|_{2}\right).
\end{align*}
Collecting the previous estimates we find that
$$
\sum_{j=1}^8I_j\leq \frac{1}{2}|\xi^\kappa,_{1111}|^2_{L^2}+C(\mathscr{E}(t))^2+C\mathscr{E}(t)+ C|\xi^\kappa,_{1111}|_{L^2}(|\xi^\kappa,_{1111}|_{L^2}+|h^\kappa,_{1111}|_{L^2})(|\xi^{\kappa}|_{2}+|h^{\kappa}|_{2}),
$$
and then,
\begin{align}\label{estimatexi}
\frac{d}{dt}|\xi,_{111}|_{L^2}^2&\leq C(1+\mathscr{E}(t))^2 -|\xi^\kappa,_{1111}|^2_{L^2}+ C|\xi^\kappa,_{1111}|_{L^2}(|\xi^\kappa,_{1111}|_{L^2}+|h^\kappa,_{1111}|_{L^2})(|\xi^{\kappa}|_{2}+|h^{\kappa}|_{2}).
\end{align}
We observe that the higher order terms on the right hand side of the previous inequality can be absorbed by the parabolic smoothing as long as $|h|_{2}+|\xi|_2$ is small enough.

\subsubsection*{The uniform time of existence}
Summing \eqref{estimatexi} and \eqref{estimateh}, using the smallness of $|h|_{2}+|\xi|_2$ and integrating we find the following inequality
$$
|\xi|_{3}^2+|h|_{3}^2\leq |\xi^\delta_0|_{3}^2+|h^\delta_0|_{3}^2 +t(1+\mathscr{E}(t))^2.
$$
Integrating \eqref{estimatephi1} and \eqref{estimatephi2} and adding them, we conclude the desired bound
$$
\mathscr{E}(t)+\frac{1}{8}\int_0^t|h^\kappa,_{1111}|^2+|\xi^\kappa,_{1111}|^2ds\leq |\xi_0|_{3}^2+|h_0|_{3}^2+tC(1+\mathscr{E}(t))^2.
$$
This bound provides us with a uniform-in-$\kappa$ time of existence $T_*$. Furthermore, in this time of existence we have that
$$
\mathscr{E}(t)\leq 2(|\xi_0|_{3}^2+|h_0|_{3}^2),
$$
and
$$
\max\{|\xi|_{2},|h|_{2}\}\leq \gamma.
$$

\subsection{Passing to the limit and uniqueness}
We have obtained uniform bounds
$$
h\in L^\infty(0,T_*;H^3(\bT)),
$$
$$
\xi\in L^\infty(0,T_*;H^3(\bT)),
$$
$$
h^\kappa\in L^2(0,T_*;H^4(\bT)),
$$
$$
\xi^\kappa\in L^2(0,T_*;H^4(\bT)),
$$
$$
\nabla \phi\in L^\infty(0,T_*;H^{2.5}(\Omega))\cap L^2(0,T_*;H^{3.5}(\Omega)),
$$
$$
\nabla \delta\psi\in L^\infty(0,T_*;H^{2.5}(\Omega))\cap L^2(0,T_*;H^{3.5}(\Omega)).
$$
With this regularity we can pass to the limit in $\kappa$ (as in \cite{majda2002vorticity}) and conclude the existence of a solution to 
\begin{align*}
\sA^\ell_j(\sA^k_j\phi,_k),_\ell&=0&&\text{ in }\Omega\times[0,T]\,,\\
\phi  &= \xi \qquad &&\text{ on }\Gamma\times[0,T]\,,\\
\lim_{x_2\rightarrow-\infty}\nabla \phi&= 0 \qquad &&\text{on }\bT\times[0,T]\,,\\
\xi_t &=-\frac{\varepsilon}{2}\sA^\ell_j\phi,_\ell\sA^k_j\phi,_k- h-\alpha \sA^\ell_2(\sA^k_2\phi,_k),_\ell\nonumber\\
&\quad+\varepsilon \sA^k_2\phi,_k\left(\sA^k_j\phi,_k\tilde{n}_j+\alpha h,_{11}\right)&&\text{ on }\Gamma\times[0,T]\,,\\
h_t&=\sA^k_j\phi,_k\tilde{n}_j+\alpha h,_{11}&&\text{ on }\Gamma\times[0,T]\,,\\
\xi &=\xi^\delta_0&&\text{ on }\bT\times\{0\}\,,\\
h&=h^\delta_0&&\text{ on }\bT\times\{0\}\,,\\
\Delta \delta\psi&= 0  &&\text{ in }\Omega\times[0,T]\,,\\
\delta\psi&= h &&\text{ on }\Gamma\times[0,T]\,,\\
\psi&= e+\delta\psi &&\text{ in }\Omega\times[0,T]\,.\end{align*}
As everything is independent of $\delta$, we can repeat the previous estimates and pass to the limit in $\delta$ in order to conclude the existence of solutions enjoying the following regularity
$$
h\in L^\infty(0,T_*;H^3(\bT))\cap L^2(0,T_*;H^4(\bT)),
$$
$$
\xi\in L^\infty(0,T_*;H^3(\bT))\cap L^2(0,T_*;H^4(\bT)),
$$
$$
\nabla \phi\in L^\infty(0,T_*;H^{2.5}(\Omega))\cap L^2(0,T_*;H^{3.5}(\Omega)),
$$
$$
\nabla \delta\psi\in L^\infty(0,T_*;H^{2.5}(\Omega))\cap L^2(0,T_*;H^{3.5}(\Omega)).
$$
The continuity in time in the highest norm is obtained in a standard way due to the parabolic effect present for $h$ and $\xi$.

Due to the regularity of the solution, the uniqueness follows the standard contradiction argument and thus, for the sake of brevity, we omit it.

\section{Proof of Theorem \ref{theorem2}: Global well-posedness in Wiener spaces}
In this section we provide appropriate estimates for the solution in Wiener spaces. Furthermore, we prove that the solution becomes analytic for positive times. In what follow we fix $\varepsilon=1$ to simplify the exposition.

\subsection{Estimates for the velocity potential} 
The first ingredient in our proof of the appropriate energy estimates in Wiener and Wiener-Sobolev spaces is a set of sharp enough estimates for the velocity potential. In particular, our goal is to prove the following estimates for $r,\lambda ,s \geq 0$
\begin{equation}\label{estPHI}
\norm{\Lambda^r \nabla \phi_2}_{\cA^{s, 1}_{\lambda}} \leqslant C(r,s) (\av{h\pare{t}}_{1, \lambda} + \av{\xi\pare{t}}_{1, \lambda})\bra{  \av{\Lambda^{r+1} \xi(t) }_{s, \lambda} + \av{ \Lambda^{r+1} h(t)}_{s, \lambda}  }  ,
\end{equation}
\begin{equation}\label{estPHI2}
\norm{\nabla \phi_2}_{\cA^{s, 2}_\lambda} \leqslant C(s)(\av{h\pare{t}}_{1, \lambda} + \av{\xi\pare{t}}_{1, \lambda})\bra{ \av{ \Lambda^2 h (t)}_{s, \lambda} +   \av{\Lambda^{2} \xi (t)}_{s, \lambda}  },
\end{equation}
\begin{equation}
\label{eq:elliptic_est_phi1}
\begin{aligned}
\av{\Lambda^r \partial_2^j \partial_i \phi_1}_{s, \lambda} & \leqslant \norm{\Lambda^r \partial_2^j \partial_i \phi_1}_{\cA^{s, 1}_\lambda} \leqslant \av{\Lambda^{r+j+1} \xi}_{s, \lambda}. 
\end{aligned}
\end{equation}
These estimates require the suitable smallness for $h$ and $\xi$. In what follows we assume that the following smallness condition holds
\begin{equation}\label{eq:smllnss1}
\av{h \pare{t}}_{1, \lambda}+ \av{\xi\pare{t}}_{1, \lambda} \leqslant C(\alpha,\varepsilon), 
\end{equation}
In the subsequent sections we will prove that, if the initial data is small enough is certain spaces, the previous smallness condition is automatically satisfied. Equipped with these estimates for the velocity potential we can, later on, establish the \emph{bona fide} estimates for the interface and the trace of the velocity potential and, as a consequence, we can conclude the estimates leading to the global well-posedness in Wiener spaces.

We recall the decomposition \eqref{phidecomposition} for the velocity potential. Then, applying Lemmas \ref{lem:product_rule_Wiener_strip} and \ref{lem:elliptic_estimate_via_symbol}, we find that, for any $ r, s, \lambda \geq 0 $, $ j\in\bN $ and $ i=1, 2 $,
\begin{equation*}
\begin{aligned}
\av{\Lambda^r \partial_2^j \partial_i \phi_1}_{s, \lambda} & \leqslant \norm{\Lambda^r \partial_2^j \partial_i \phi_1}_{\cA^{s, 1}_\lambda} \leqslant \av{\Lambda^{r+j+1} \xi}_{s, \lambda}. 
\end{aligned}
\end{equation*}

Now we have to find appropriate estimates for $ \phi_2 $ solving the equation \eqref{eq:eq_phi2}. This equation contains the nonlinear contributions and, as a consequence, it is more challenging to estimate.

We start invoking Proposition \ref{prop:elliptic_estimates} and applying it to \eqref{eq:eq_phi2} in order to ensure that
$$
\norm{\nabla \phi_2}_{\cA^{0, 1}_{\lambda}}\leq 12\norm{(J \sA \sA^\intercal -I)(\nabla \phi_1+\nabla \phi_2)}_{\cA^{0, 1}_{\lambda}}.
$$
Using Lemma \ref{lem:product_rule_Wiener_strip}, we find that
\begin{align*}
\norm{\nabla \phi_2}_{\cA^{0, 1}_{\lambda}}&\leq 24\norm{J \sA \sA^\intercal -I}_{\cA^{0, 1}_{\lambda}}\norm{\nabla \phi_1}_{\cA^{0, 1}_{\lambda}}+24\norm{J \sA \sA^\intercal -I}_{\cA^{0, 1}_{\lambda}}\norm{\nabla \phi_2}_{\cA^{0, 1}_{\lambda}}\\
&\leq 24\norm{J \sA \sA^\intercal -I}_{\cA^{0, 1}_{\lambda}}\norm{\nabla \phi_1}_{\cA^{0, 1}_{\lambda}}+240 \av{ \Lambda h}_{0, \lambda}\norm{\nabla \phi_2}_{\cA^{0, 1}_{\lambda}}.
\end{align*}
Thus,
\begin{align}\label{termaux}
\norm{\nabla \phi_2}_{\cA^{0, 1}_{\lambda}}\leq \frac{24\norm{J \sA \sA^\intercal -I}_{\cA^{0, 1}_{\lambda}}\norm{\nabla \phi_1}_{\cA^{0, 1}_{\lambda}}}{(1-240 \av{ \Lambda h}_{0, \lambda})}.
\end{align}

We apply Proposition in \ref{prop:elliptic_estimates} to \eqref{eq:eq_phi2} obtaining the estimate
\begin{equation}\label{eq:bound_elliptic_est_phi2}
\norm{\Lambda^r \nabla\phi_2}_{\cA^{s, 1}_{\lambda}}\leqslant 12 \norm{\Lambda^r \bra{ \pare{J \sA \sA^\intercal - I}\nabla\pare{\phi_1+\phi_2}}}_{\cA^{s, 1}_{\lambda}}, 
\end{equation}
where $ s\geq 1 $ and $ r, \lambda \geq 0 $. 

Let us now apply the product rule stated in Lemma \ref{lem:product_rule_Wiener_strip} to deduce
\begin{multline}\label{eq:application_elliptic_estimates}
\norm{\Lambda^r \nabla \phi_2}_{\cA^{s, 1}_{\lambda}}
%----------------------------------------------------------------
 \leqslant 24 \ \ck_r \ck_s \bigg[ \norm{J \sA \sA^\intercal - I}_{\cA^{0, 1}_{\lambda}} \pare{\norm{\Lambda^r \nabla \phi_1}_{\cA^{s, 1}_{\lambda}} + \norm{\Lambda^r \nabla \phi_2}_{\cA^{s, 1}_{\lambda}}} \\
 %----------------------------------------------------------------
  + \norm{\Lambda^r \pare{ J \sA \sA^\intercal - I }}_{\cA^{s, 1}_{\lambda}} \pare{\norm{\nabla \phi_1}_{\cA^{0, 1}_{\lambda}} + \norm{\nabla \phi_2}_{\cA^{0, 1}_{\lambda}}} \bigg]. 
\end{multline}

At this point it is easy to find (see also \cite{GGS19} for more details) that if $ \av{h}_{1, \lambda} < \pare{16 \ck_r \ck_s \varepsilon }^{-1} $ then the following bound holds true
\begin{equation}\label{eq:est_Q}
\norm{\Lambda^r \pare{ J \sA \sA^\intercal -I }}_{\cA^{s, 1}_{\lambda}} \leqslant 10 \av{ \Lambda^{r+1} h}_{s, \lambda}.
\end{equation}

We are finally now able to bound the right hand side of \eqref{eq:bound_elliptic_est_phi2}:

Using \eqref{eq:est_Q} in order to deduce that
\begin{equation*}
24 \ \ck_r \ck_s \norm{J \sA \sA^\intercal -I}_{\cA^{0, 1}_{\lambda}} \norm{\Lambda^r \nabla \phi_2}_{\cA^{s, 1}_{\lambda}} \leqslant 240  \ \ck_r \ck_s    \av{h}_{1, \lambda} \norm{\Lambda^r \nabla \phi_2}_{\cA^{s, 1}_{\lambda}}. 
\end{equation*}
We have that \eqref{eq:bound_elliptic_est_phi2} together with the previous inequality leads to
\begin{multline}\label{eq:bound_elliptic_est_phi2_2}
\pare{ 1-240 \ \ck_r \ck_s 
    \av{h}_{1, \lambda}
}
\norm{\Lambda^r \nabla \phi_2}_{\cA^{s, 1}_{\lambda}}\\
%----------------------------------------------------------------
 \leqslant 24 \ \ck_r \ck_s \bra{ \norm{J \sA \sA^\intercal -I}_{\cA^{0, 1}_{\lambda}} \norm{\Lambda^r \nabla \phi_1}_{\cA^{s, 1}_{\lambda}}  + \norm{ \Lambda^r \pare{ J \sA \sA^\intercal -I }}_{\cA^{s, 1}_{\lambda}} \norm{\nabla \phi_1}_{\cA^{0, 1}_{\lambda}} } \\
 + \norm{\Lambda^r \pare{ J \sA \sA^\intercal - I }}_{\cA^{s, 1}_{\lambda}} \norm{\nabla \phi_2}_{\cA^{0, 1}_{\lambda}} \bigg]. 
\end{multline}
We apply now \eqref{eq:elliptic_est_phi1}, \eqref{termaux} and \eqref{eq:est_Q} to \eqref{eq:bound_elliptic_est_phi2_2} in order to obtain an inequality which involves norms of $ \phi_2 $ and $ h $ only
\begin{align*}
\pare{1-
240  \ \ck_r \ck_s  \av{h}_{1, \lambda}
}
\norm{\Lambda^r \nabla \phi_2}_{\cA^{s, 1}_{\lambda}}
 %----------------------------------------------------------------
&\leqslant
  C(r,s) \bigg{[}     \av{h}_{1, \lambda} \av{\Lambda^{r+1} \xi }_{s, \lambda} + \av{ \Lambda^{r+1} h}_{s, \lambda} \av{\Lambda \xi }_{0, \lambda} 
\\
&\quad +  \frac{ \av{\Lambda \xi}_{0, \lambda} \av{\Lambda^{r+1} h}_{s, \lambda} \av{\Lambda h}_{0, \lambda} }{1-240  \av{h}_{1, \lambda}} \bigg{]},
\end{align*}
where $C(r,s)$ is a (computable) constant that only depends on $r$ and $s$ and may change from line to line. Using hence the smallness of $h$, we prove that
\begin{equation}\label{eq:elliptic_est_phi2_trace}
\norm{\Lambda^r \nabla \phi_2}_{\cA^{s, 1}_{\lambda}} \leqslant C(r,s) \bra{     \av{h}_{1, \lambda} \av{\Lambda^{r+1} \xi }_{s, \lambda} + \av{ \Lambda^{r+1} h}_{s, \lambda} \av{ \xi }_{1, \lambda} }  .
\end{equation}
The estimate \eqref{eq:elliptic_est_phi2_trace} allow us to control all terms appearing in \eqref{eq:aux1} with the exception of $ \partial_2^2\phi_2 $. It is in this context that we use the higher order elliptic estimates proved in the Proposition in \ref{prop:elliptic_estimates}. Thus applying Proposition in \ref{prop:elliptic_estimates} we recover the inequality
\begin{align}\label{eq:est_nablaphi2_As2lambda_1}
\norm{\nabla \phi_2}_{\cA^{s, 2}_\lambda}& \leqslant 12 \norm{\Lambda \bra{(J \sA \sA^\intercal -I)\pare{\phi_1+\phi_2}}}_{\cA^{s, 1}_\lambda} +  4 \norm{(J \sA \sA^\intercal -I)\pare{\phi_1+\phi_2}}_{\cA^{s, 2}_\lambda}\nonumber\\
& = I_1^{s, \lambda} + I_2^{s, \lambda},  
\end{align} 
which holds for $ s, \lambda\geq 0 $.

The term $ I_1^{s, \lambda} $ is easier to control. To estimate this term we apply Lemma \ref{lem:product_rule_Wiener_strip} with $ r=1 $ in order to deduce
\begin{align}
I_1^{s, \lambda} \leqslant 24 \ck_s \bigg{[} \norm{\Lambda (J \sA \sA^\intercal -I)}_{\cA^{s, 1}_\lambda}\pare{ \norm{\nabla\phi_1}_{\cA^{0, 1}_\lambda} + \norm{\nabla\phi_2}_{\cA^{0, 1}_\lambda}  }
\nonumber\\+
 \norm{ (J \sA \sA^\intercal -I)}_{\cA^{0, 1}_\lambda}\pare{ \norm{\Lambda \nabla\phi_1}_{\cA^{s, 1}_\lambda} + \norm{\Lambda\nabla\phi_2}_{\cA^{s, 1}_\lambda}  }
\bigg{]}\label{eq:est_Islambda1_1}, 
\end{align}
Recalling \eqref{eq:est_Q}, we obtain that
\begin{equation}\label{eq:est_Islambda1_2}
\norm{(J \sA \sA^\intercal -I)}_{\cA^{0, 1}_\lambda}  \leqslant 10 \av{ h}_{1, \lambda}, \quad
\norm{\Lambda (J \sA \sA^\intercal -I)}_{\cA^{s, 1}_\lambda}  \leqslant 10 \av{\Lambda^2 h}_{s, \lambda}. 
\end{equation}

Similarly, using the formula of $\phi_1$ and \eqref{eq:elliptic_est_phi1}, we find that
\begin{equation}\label{eq:est_Islambda1_3}
\norm{\nabla \phi_1}_{\cA^{0, 1}_\lambda}  \leqslant 2 \av{\xi}_{1, \lambda}, \quad
\norm{\Lambda \nabla \phi_1}_{\cA^{s, 1}_\lambda}  \leqslant 2 \av{\Lambda^2 \xi}_{s, \lambda}. 
\end{equation}

We use the smallness hypothesis together with \eqref{termaux}, \eqref{eq:elliptic_est_phi2_trace}, \eqref{eq:est_Islambda1_2}, \eqref{eq:est_Islambda1_2} in \eqref{eq:est_Islambda1_1} in order to obtain
\begin{equation*}
I_1^{s, \lambda}\leqslant C(s,\varepsilon) \bra{ \av{ \Lambda^2 h }_{s, \lambda} \pare{1+ \av{h}_{1, \lambda}} \av{ \xi }_{1, \lambda} + \av{h}_{1, \lambda} \pare{ \pare{1+ \av{h}_{1, \lambda}  }\av{\Lambda^{2} \xi }_{s, \lambda} + 2 \av{ \xi }_{1, \lambda} \av{ \Lambda^{2} h}_{s, \lambda} }}   ,  
\end{equation*}
which after algebraic manipulation gives us
\begin{equation}\label{lem:Islambda1}
I_1^{s, \lambda} \leqslant C(s,\varepsilon)(\av{h\pare{t}}_{1, \lambda} + \av{\xi\pare{t}}_{1, \lambda})\bra{ \av{ \Lambda^2 h }_{s, \lambda} +   \av{\Lambda^{2} \xi }_{s, \lambda}  } .
\end{equation}

Before starting to prove a suitable bound for the term $ I_2^{s, \lambda} $ we need approriate estimates for
\begin{equation*}
\norm{\Lambda^r (J \sA \sA^\intercal -I)}_{\cA^{s, 2}_\lambda}. 
\end{equation*}
First, let us remark that
\begin{equation*}
\norm{\Lambda^r (J \sA \sA^\intercal -I)}_{\cA^{s, 2}_\lambda} = \norm{\Lambda^r \partial_2  (J \sA \sA^\intercal -I)}_{\cA^{s, 1}_\lambda}. 
\end{equation*}
Now, we compute 
% \begin{equation*}
% \partial_2 (J \sA \sA^\intercal -I) = \pare{
%\begin{array}{cc}
%\partial_2^2 \delta\psi & -\partial_{12} \delta\psi \\[2mm]
%-\partial_{12} \delta\psi & \displaystyle \frac{-\partial_2^2 \delta\psi + 2\varepsilon\partial_1 \delta\psi\partial_{12} \delta\psi}{1+\varepsilon \partial_2 \delta\psi}
%-\frac{\varepsilon\pare{ -\partial_2 \delta\psi +\varepsilon \pare{\partial_1 \delta\psi}^2 }\partial_2^2 \delta\psi }{\pare{ 1+\varepsilon \partial_2 \delta\psi }^2 }
%\end{array}
%},
% \end{equation*}
 \begin{equation*}
 \partial_2 (J \sA \sA^\intercal -I) = \pare{
\begin{array}{cc}
\partial_2^2 \delta\psi & -\partial_{12} \delta\psi \\[2mm]
-\partial_{12} \delta\psi & \displaystyle \frac{-\partial_2^2 \delta\psi + 2\partial_1 \delta\psi\partial_{12} \delta\psi}{1+ \partial_2 \delta\psi}
-\frac{\pare{ -\partial_2 \delta\psi + \pare{\partial_1 \delta\psi}^2 }\partial_2^2 \delta\psi }{\pare{ 1+ \partial_2 \delta\psi }^2 }
\end{array}
},
 \end{equation*}

and we use the identities 
$$ 
\pare{1+x}^{-1} =  \sum_{n=0}^\infty \pare{-x}^n  $$ 
and 
$$ 
\pare{1+x}^{-2} = \sum_{n=0}^\infty \pare{-1}^n \pare{n+1} x^n,
$$
both valid for $  \av{x}<1 $ in order to express it as 
$$
\partial_2 (J \sA \sA^\intercal -I) = \sum_{n=0}^\infty \tilde{Q}_{\pare{n}},
$$
with
\begin{align*}
\tilde{Q}_{\pare{0}} & = \pare{
\begin{array}{cc}
\partial_2^2 \delta\psi & -\partial_{12} \delta\psi \\[2mm]
-\partial_{12} \delta\psi & -\partial_2^2 \delta\psi + 2\partial_1 \delta\psi\partial_{12} \delta\psi
-\pare{ -\partial_2 \delta\psi +\pare{\partial_1 \delta\psi}^2 }\partial_2^2 \delta\psi 
\end{array}
}, \\
%------------------------------------------------------
\tilde{Q}_{\pare{n}} & = \pare{
\begin{array}{cc}
0 & 0  \\
0 & \pare{ -\partial_2^2 \delta\psi + 2\partial_1 \delta\psi\partial_{12} \delta\psi }\pare{-\partial_2\delta\psi}^n
-\pare{ n + 1 } \pare{ \pare{ -\partial_2 \delta\psi +\pare{\partial_1 \delta\psi}^2 }\partial_2^2 \delta\psi  }\pare{-\partial_2\delta\psi}^n
\end{array}
}, && n\geq 1.
\end{align*}
We can apply Lemma \ref{lem:product_rule_Wiener_strip} together with the inequality
\begin{equation*}
 \norm{\Lambda^r \partial_1^k\partial_2^l \delta\psi}_{\cA^{s, 1}_\lambda}\leqslant\av{\Lambda^{r+k+l} h}_{s, \lambda}, 
\end{equation*}
which stems from the explicit definition of $ \delta\psi $ and Lemma \ref{lem:elliptic_estimate_via_symbol}. Applying the product rules of \ref{lem:product_rule_Wiener_strip} it is possible to produce the bound
\begin{align}
&\norm{\Lambda^r \pare{ \pare{-\partial_2\delta\psi}^n\pare{- \partial_2^2 \delta\psi + 2\partial_1 \delta\psi \ \partial_{12}\delta\psi}}}_{\cA^{s, 1}_\lambda } \nonumber\\
 \leqslant&  C(r,s) \Big\{
%--------------------------------------------------------------
\norm{\Lambda^r\pare{ \pare{\partial_2 \delta\psi}^n }}_{\cA^{s, 1}_\lambda} \bra{ \norm{\partial_2^2\delta\psi}_{\cA^{0, 1}_\lambda} + \norm{\partial_1 \delta\psi}_{\cA^{0, 1}_\lambda} \norm{\partial_{12} \delta\psi}_{\cA^{0, 1}_\lambda} } \nonumber\\
%--------------------------------------------------------------
& + \norm{ \pare{\partial_2 \delta\psi}^n }_{\cA^{0, 1}_\lambda} \bra{ \norm{\Lambda^r \partial_2^2\delta\psi}_{\cA^{s, 1}_\lambda} + \pare{ \norm{\Lambda^r \partial_1 \delta\psi}_{\cA^{s, 1}_\lambda} \norm{\partial_{12} \delta\psi}_{\cA^{0, 1}_\lambda} +\norm{\partial_1 \delta\psi}_{\cA^{0, 1}_\lambda} \norm{\Lambda^r \partial_{12} \delta\psi}_{\cA^{s, 1}_\lambda}  }} \Big\}, \nonumber\\
%--------------------------------------------------------------
%--------------------------------------------------------------
 \leqslant&  \ C(r,s)^{n+1} \Big\{
%--------------------------------------------------------------
\norm{\partial_2 \delta\psi}_{\cA^{0, 1}_\lambda}^{n-1} \norm{\Lambda^r\partial_2 \delta\psi}_{\cA^{s, 1}_\lambda} \bra{ \norm{\partial_2^2\delta\psi}_{\cA^{0, 1}_\lambda} + \norm{\partial_1 \delta\psi}_{\cA^{0, 1}_\lambda} \norm{\partial_{12} \delta\psi}_{\cA^{0, 1}_\lambda} } \nonumber\\
%--------------------------------------------------------------
& + \norm{ \partial_2 \delta\psi }_{\cA^{0, 1}_\lambda}^n \bra{ \norm{\Lambda^r \partial_2^2\delta\psi}_{\cA^{s, 1}_\lambda} + \pare{ \norm{\Lambda^r \partial_1 \delta\psi}_{\cA^{s, 1}_\lambda} \norm{\partial_{12} \delta\psi}_{\cA^{0, 1}_\lambda} +\norm{\partial_1 \delta\psi}_{\cA^{0, 1}_\lambda} \norm{\Lambda^r \partial_{12} \delta\psi}_{\cA^{s, 1}_\lambda}  }} \Big\}, \nonumber\\
%--------------------------------------------------------------
%--------------------------------------------------------------
 \leqslant&  \ C(r,s)^{n+1}  \Big\{
%--------------------------------------------------------------
\av{\Lambda h}_{0, \lambda}^{n-1} \av{\Lambda^{1+r} h}_{s, \lambda} \bra{ \av{\Lambda^2 h}_{0, \lambda}+ \av{\Lambda h}_{0, \lambda} \av{\Lambda^2 h}_{0, \lambda}} \nonumber\\
%--------------------------------------------------------------
& + \av{\Lambda h}_{0, \lambda}^n \bra{ \av{\Lambda^{2+r} h}_{s, \lambda} + \pare{ \av{\Lambda^{1+r}h}_{s, \lambda} \av{\Lambda^2 h}_{0, \lambda} +\av{\Lambda h}_{0, \lambda} \av{\Lambda^{2+r} h }_{s, \lambda} }} \Big\}, \label{eq:Qtilde_1}
\end{align}
since $ h $ preserve its (zero) average we interpolate and use Poincar\'e inequality in order to deduce the bound
\begin{align}\label{eq:Qtilde_2}
\av{\Lambda^{1+r}h}_{s, \lambda}\av{\Lambda^2 h}_{0, \lambda}\leqslant \av{\Lambda^{2+r} h}_{s, \lambda}\av{\Lambda h}_{0, \lambda}, && s\geq 0, \ r\geq 1. 
\end{align}
This estimate applied to \eqref{eq:Qtilde_1} gives
\begin{equation}\label{eq:Qtilde_2.1}
\begin{aligned}
&\norm{\Lambda^r \pare{ \pare{-\partial_2\delta\psi}^n\pare{\partial_2^2 \delta\psi + 2\partial_1 \delta\psi \ \partial_{12}\delta\psi}}}_{\cA^{s, 1}_\lambda } \\
%------------------------------------------------------------------
 \leqslant&  \ C(r,s)^{n+1}\ \av{h}_{1, \lambda}^{n} \Big\{
%--------------------------------------------------------------
  \pare{  1 + \av{\Lambda h}_{0, \lambda} }    +   \pare{ 1 +  \av{\Lambda h}_{0, \lambda} } \Big\} \av{\Lambda^{2+r} h}_{s, \lambda}, \\
%------------------------------------------------------------------
\leqslant&  \ C(r,s)^{n+1}\ \av{h}_{1, \lambda}^{n} \pare{ 1 + \av{ h}_{1 , \lambda} } \av{\Lambda^{2+r} h}_{s, \lambda} .
\end{aligned}
\end{equation}

We perform now similar computations as the ones carried out in \eqref{eq:Qtilde_1} in order to deduce the following bound
\begin{equation}\label{eq:Qtilde_3}
\norm{\Lambda^r \pare{\pare{-\partial_2\delta\psi}^n\pare{ \pare{ -\partial_2 \delta\psi +\pare{\partial_1 \delta\psi}^2 }\partial_2^2 \delta\psi  }}}_{\cA^{s, 1}_\lambda} 
\leqslant C(r,s)^{n+1} \av{h}_{1, \lambda}^{n+1} \bra{1+\av{h}_{1, \lambda}}\av{\Lambda^{2+r} h }_{s, \lambda}. 
\end{equation}

At this point we can combine \eqref{eq:Qtilde_2.1} and \eqref{eq:Qtilde_3} in order to deduce that if 
\begin{equation*}
\av{h}_{1, \lambda}\leqslant C(r,s), 
\end{equation*}
then
\begin{equation*}
\begin{aligned}
\norm{\Lambda^r \tilde{Q}_{\pare{0}}}_{\cA^{s, 1}_\lambda} & \leqslant 4 \av{\Lambda^{2+r} h}_{s, \lambda}, \\
%---------------------------------------------
\norm{\Lambda^r \tilde{Q}_{\pare{n}}}_{\cA^{s, 1}_\lambda} & \leqslant \av{\Lambda^{2+r} h}_{s, \lambda},  & n \geq 1.
\end{aligned}
\end{equation*}
Since 
\begin{equation*}
\norm{\Lambda^r (J \sA \sA^\intercal -I)}_{\cA^{s, 2}_\lambda} = \norm{\Lambda^r  \partial_2 (J \sA \sA^\intercal -I)}_{\cA^{s, 1}_\lambda}  \leqslant \sum_n  \norm{\Lambda^r \tilde{Q}_{\pare{n}}}_{\cA^{s, 1}_\lambda}  \leqslant 5 \av{\Lambda^{2+r}h}_{s, \lambda}, 
\end{equation*}
we conclude the following estimate
\begin{equation*}
\norm{\Lambda^r (J \sA \sA^\intercal -I)}_{\cA^{s, 2}_\lambda}  \leqslant 5 \av{\Lambda^{2+r}h}_{s, \lambda}. 
\end{equation*}

We focus now on the term $ I_2^{s, \lambda} $. We apply Lemma \ref{lem:product_rule_Wiener_strip}  in order to obtain
\begin{multline}\label{eq:I2slambda_1}
I_2^{s, \lambda} \leqslant C(s) \bigg[ 
\norm{(J \sA \sA^\intercal -I)}_{\cA^{s, 2}_\lambda} \pare{\norm{\nabla\phi_1}_{\cA^{0, 1}_\lambda} + \norm{\nabla\phi_2}_{\cA^{0, 1}_\lambda}} \\  
+ \norm{(J \sA \sA^\intercal -I)}_{\cA^{0, 1}_\lambda} \pare{\norm{\nabla\phi_1}_{\cA^{s, 2}_\lambda} + \norm{\nabla\phi_2}_{\cA^{s, 2}_\lambda}} \\
%----------------------------------------------------------------
+ \norm{(J \sA \sA^\intercal -I)}_{\cA^{s, 1}_\lambda} \pare{\norm{\nabla\phi_1}_{\cA^{0, 2}_\lambda} + \norm{\nabla\phi_2}_{\cA^{0, 2}_\lambda}}\\
+\norm{(J \sA \sA^\intercal -I)}_{\cA^{0, 2}_\lambda} \pare{\norm{\nabla\phi_1}_{\cA^{s, 1}_\lambda} + \norm{\nabla\phi_2}_{\cA^{s, 1}_\lambda}} \bigg]. 
\end{multline}
The term $ \norm{\nabla \phi_2 }_{\cA^{0, 2}_\lambda} $ is a high order term taht we have to contro in terms of well-behaved quantities that we know how to bound. Let us hence apply again Proposition in \ref{prop:elliptic_estimates} in order to obtain
\begin{equation*}
\norm{\nabla \phi_2}_{\cA^{0, 2}_\lambda} \leqslant 12 \norm{\Lambda\bra{(J \sA \sA^\intercal -I)\nabla\pare{\phi_1 + \phi_2}}}_{\cA^{0, 1}_\lambda} + 4 \norm{(J \sA \sA^\intercal -I)\nabla\pare{\phi_1 + \phi_2}}_{\cA^{0, 2}_\lambda}, 
\end{equation*}
and we use the product rules of \ref{lem:product_rule_Wiener_strip} in order to deduce the bound 
\begin{multline*}
\norm{(J \sA \sA^\intercal -I)\nabla\pare{\phi_1 + \phi_2}}_{\cA^{0, 2}_\lambda}\leq   \  C \bigg{[} \norm{(J \sA \sA^\intercal -I)}_{\cA^{0, 2}_\lambda} \pare{ \norm{\nabla\phi_1}_{\cA^{0, 1}_\lambda} + \norm{\nabla\phi_2}_{\cA^{0, 1}_\lambda}}
\\
\quad+
\norm{(J \sA \sA^\intercal -I)}_{\cA^{0, 1}_\lambda} \pare{ \norm{\nabla\phi_1}_{\cA^{0, 2}_\lambda} + \norm{\nabla\phi_2}_{\cA^{0, 2}_\lambda}}\bigg{]}
, 
\end{multline*}
and 
\begin{multline*}
\norm{\Lambda\bra{(J \sA \sA^\intercal -I)\nabla\pare{\phi_1 + \phi_2}}}_{\cA^{0, 1}_\lambda} \leqslant C  \left[ \norm{\Lambda (J \sA \sA^\intercal -I)}_{\cA^{0, 1}_\lambda}\pare{ \norm{\nabla\phi_1}_{\cA^{0, 1}_\lambda} + \norm{\nabla\phi_2}_{\cA^{0, 1}_\lambda} } \right. \\
%--------------------------------------------------
+ \left.  \norm{ (J \sA \sA^\intercal -I)}_{\cA^{0, 1}_\lambda}\pare{ \norm{\Lambda \nabla\phi_1}_{\cA^{0, 1}_\lambda} + \norm{\Lambda \nabla\phi_2}_{\cA^{0, 1}_\lambda} } \right] , 
\end{multline*}
from which we deduce that
\begin{multline}
\label{eq:nablaphi2_02_1}
\norm{\nabla \phi_2}_{\cA^{0, 2}_\lambda} \leqslant \frac{C}{1-\norm{(J \sA \sA^\intercal -I)}_{\cA^{0, 1}_\lambda}}\\
%-------------------------------------------------------
\bigg\{ \norm{\Lambda (J \sA \sA^\intercal -I)}_{\cA^{0, 1}_\lambda}\pare{ \norm{\nabla\phi_1}_{\cA^{0, 1}_\lambda} + \norm{\nabla\phi_2}_{\cA^{0, 1}_\lambda} }
+  \norm{ (J \sA \sA^\intercal -I)}_{\cA^{0, 1}_\lambda}\pare{ \norm{\Lambda \nabla\phi_1}_{\cA^{0, 1}_\lambda} + \norm{\Lambda \nabla\phi_2}_{\cA^{0, 1}_\lambda} }\\
%-------------------------------------------------------
 + 2 \bra{ \norm{(J \sA \sA^\intercal -I)}_{\cA^{0, 2}_\lambda} \pare{ \norm{\nabla\phi_1}_{\cA^{0, 1}_\lambda} + \norm{\nabla\phi_2}_{\cA^{0, 1}_\lambda}}
+
\norm{(J \sA \sA^\intercal -I)}_{\cA^{0, 1}_\lambda}  \norm{\nabla\phi_1}_{\cA^{0, 2}_\lambda} 
}\bigg\}. 
\end{multline}
Using the smallness hypothesis \eqref{eq:smllnss1} together with the estimates \eqref{eq:elliptic_est_phi1} and \eqref{eq:est_Q}, we obtain the bounds
\begin{equation}
\label{eq:nablaphi2_02_2}
\begin{aligned}
\norm{\Lambda^r (J \sA \sA^\intercal -I)}_{\cA^{s, 1}_\lambda} & \leqslant 10 \av{\Lambda^{1+r}h}_{s, \lambda}, \\
%---------------------------------------
\norm{\Lambda^r\nabla \phi_1}_{\cA^{s, j}_\lambda} & \leqslant 2 \av{\Lambda^{j+r}\xi}_{s, \lambda}, \\
%---------------------------------------
\norm{\Lambda^r \nabla \phi_2}_{\cA^{s, 1}_\lambda} & \leqslant C\pare{\av{h}_{1, \lambda} + \av{\xi}_{1, \lambda}} \pare{\av{\Lambda^{1+r}\xi}_{s, \lambda} + \av{\Lambda^{1+r} h}_{s, \lambda}}, 
\end{aligned}
\end{equation}
where $ r, s, \lambda \geq 0 $ and $ j\in\bN $. We use now \eqref{eq:nablaphi2_02_2} in \eqref{eq:nablaphi2_02_1} and after suitable simplification we find the estimate
\begin{equation}
\label{eq:nablaphi2_02_3}
\begin{aligned}
\norm{\nabla\phi_2}_{\cA^{0, 2}_\lambda} & \leqslant C\pare{\av{h}_{1, \lambda} + \av{\xi}_{1, \lambda}}\pare{\av{\Lambda^2 h}_{0, \lambda} + \av{\Lambda^2 \xi }_{0, \lambda}}. %\\
%------------------------------------
%& \leqslant c_0\nu \pare{\av{\Lambda^2 h}_{0, \lambda} + \av{\Lambda^2 \xi }_{0, \lambda}} . 
\end{aligned}
\end{equation}

%Rearranging terms we can transform \eqref{eq:I2slambda_1} into
%\begin{multline}\label{eq:I2slambda_2}
%I_2^{s, \lambda} \leqslant 8 \ck_s \bigg\{ 
%\norm{Q\pare{\nabla \sigma}}_{\cA^{s, 2}_\lambda} \pare{\norm{\nabla\phi_1}_{\cA^{0, 1}_\lambda} + \norm{\nabla\phi_2}_{\cA^{0, 1}_\lambda}}   
%+ \norm{Q\pare{\nabla \sigma}}_{\cA^{0, 1}_\lambda} \norm{\nabla\phi_1}_{\cA^{s, 2}_\lambda}  \\
%%----------------------------------------------------------------
%+ \norm{Q\pare{\nabla \sigma}}_{\cA^{s, 1}_\lambda} \pare{\norm{\nabla\phi_1}_{\cA^{0, 2}_\lambda} + \norm{\nabla\phi_2}_{\cA^{0, 2}_\lambda}}
%+\norm{Q\pare{\nabla \sigma}}_{\cA^{0, 2}_\lambda} \pare{\norm{\nabla\phi_1}_{\cA^{s, 1}_\lambda} + \norm{\nabla\phi_2}_{\cA^{s, 1}_\lambda}} \bigg\}
% \\
%%----------------------------------------------------------------
%+ \norm{Q\pare{\nabla \sigma}}_{\cA^{0, 1}_\lambda} \norm{\nabla\phi_2}_{\cA^{s, 2}_\lambda}
%,
%\end{multline}
%so that we can isolate the high-order contribution deriving from $ \phi_2 $. \\

We use the inequality \eqref{eq:technical_interpolation_estimate} to obtain
\begin{equation*}
\pare{\av{\Lambda^2 \xi}_{0, \lambda} + \av{\Lambda^2 h}_{0, \lambda}} \pare{\av{\Lambda \xi}_{s, \lambda} + \av{\Lambda h}_{s, \lambda}} \leqslant 2^{1+\frac{1}{1+s}} \pare{\av{ \xi}_{1, \lambda} + \av{h}_{1, \lambda}} \pare{\av{\Lambda^2 \xi}_{s, \lambda} + \av{\Lambda^2 h}_{s, \lambda}}.
\end{equation*}

We use now the estimates \eqref{eq:smllnss1}, \eqref{eq:est_Q}, \eqref{eq:nablaphi2_02_2} and \eqref{eq:nablaphi2_02_3} together with the previous inequality to conclude that
\begin{equation*}
I^{s, \lambda}_2 \leqslant C(s)\pare{\av{ \xi}_{1, \lambda} + \av{h}_{1, \lambda}}  \pare{\av{\Lambda^2 \xi}_{s, \lambda} + \av{\Lambda^2 h}_{s, \lambda}} 
+ 10\av{h}_{1, \lambda} \norm{\nabla\phi_2}_{\cA^{s, 2}_\lambda}.
\end{equation*}
Thus, rearranging the terms in \eqref{eq:est_nablaphi2_As2lambda_1} and using the smallness hypothesis \eqref{eq:smllnss1}, we find that
$$
\norm{\nabla \phi_2}_{\cA^{s, 2}_\lambda} \leqslant C(s)(\av{h\pare{t}}_{1, \lambda} + \av{\xi\pare{t}}_{1, \lambda})\bra{ \av{ \Lambda^2 h (t)}_{s, \lambda} +   \av{\Lambda^{2} \xi (t)}_{s, \lambda}  }
$$
as desired. Furthermore, we use the trace estimates of \ref{lem:product_rule_Wiener_strip} in order to obtain
\begin{equation*}
\av{\partial_2^2\phi_2}_{s, \lambda}\leqslant \norm{\partial_2^2\phi_2}_{\cA^{s, 1}_\lambda}\leqslant \norm{\nabla \phi_2}_{\cA^{s, 2}_\lambda}\leq C(s)(\av{h\pare{t}}_{1, \lambda} + \av{\xi\pare{t}}_{1, \lambda})\bra{ \av{ \Lambda^2 h (t)}_{s, \lambda} +   \av{\Lambda^{2} \xi (t)}_{s, \lambda}  }. 
\end{equation*}
Equipped with these estimates we can now perform the analysis of the equations for $h$ and $\xi$.

\subsection{Estimates for the evolution equations}
We define
\begin{align} \label{eq:muc0}
\mu \in \left[ 0, \frac{\alpha}{2} \right). 
\end{align}
We compute that
\begin{equation*}
\begin{aligned}
\ddt \av{f\pare{t}}_{1, \mu t} & = \pare{\sum_{n=-\infty}^\infty \pare{1+\av{n}} e^{\mu t \av{n}} \partial_t \av{\hat{f}\pare{n, t}} } + \mu \pare{\sum_{n=-\infty}^\infty \pare{1+\av{n}} \av{n} e^{\mu t \av{n}}  \av{\hat{f}\pare{n, t}} }, \\
%-------------------------------------------------------------
& \leqslant   \pare{\sum_{n=-\infty}^\infty \pare{1+\av{n}} e^{\mu t \av{n}} \partial_t \av{\hat{f}\pare{n, t}} } + \mu \av{\Lambda^2 f}_{1, \mu t} . 
\end{aligned}
\end{equation*}
Now we perform an $ A^1_{\mu t} $ energy estimate on the evolution equation of $ \xi $ in \eqref{eq:ALE}. Thus, we find that

\begin{multline}\label{eq:EEpsi1}
\ddt \av{\xi}_{1, \mu t} + \left(\alpha-\mu\right) \av{\Lambda^2 \xi}_{1, \mu t} \leqslant  \av{h}_{1, \mu t} + \frac{1}{2}\av{\sA^\ell_j\phi,_\ell\sA^k_j\phi,_k}_{1, \mu t}\\
+\av{\frac{\phi,_2}{1+\Lambda h}\left(\sA^k_j\phi,_k\tilde{n}_j+ \alpha h,_{11}\right)}_{1, \mu t} \\
%-------------------------------------------
+ \alpha\av{\frac{\pare{2 + \Lambda h }\Lambda h}{\pare{1+\Lambda h}^2} \ \Lambda^2 \xi}_{1, \mu t} + \alpha\av{\frac{\Lambda^2 h}{\pare{1+\Lambda h}^3} \Lambda \xi}_{1, \mu t }  \\
+ \alpha\av{\frac{\phi_2,_{22}}{\pare{1+\Lambda h}^2}  }_{1, \mu t} + \alpha\av{\frac{\Lambda^2 h}{\pare{1+\Lambda h}^3} \phi_2,_2}_{1, \mu t}.
\end{multline}
Similarly, we find that
\begin{equation}\label{eq:EEh1}
\ddt \av{h}_{1, \mu} + \left(\alpha-\mu\right)\av{\Lambda^2 h}_{1, \mu t} = \av{\sA_1^k\partial_k \pare{ \phi_1 +\phi_2 } \ \partial_1 h}_{1, \mu t} + \av{\sA_2^2 \partial_2 \pare{ \phi_1 +\phi_2 }}_{1, \mu t}. 
\end{equation}

Before starting to compute suitable bounds for the nonlinear energy contributions in \eqref{eq:EEpsi1} and \eqref{eq:EEh1} we need the analog of estimate \eqref{eq:est_Q} in Wiener spaces. To prove such estimate, we observe that, from the explicit definition of $ \sA $ given in section \ref{sec:matA} and the definition of the function $ \mathsf{G} $ given in \eqref{eq:def_G}, we can rewrite the trace of $ A $ in $ x_2 =0 $ as
\begin{align*}
\left. \sA\right|_{x_2 =0} & = I + \widetilde{\sA} , \\
%--------------------------------------------
\widetilde{\sA} & = \pare{
\begin{array}{cc}
0 & 0 \\[3mm]
 \partial_1 h \pare{\mathsf{G}\pare{\Lambda h} - 1} & - \mathsf{G}\pare{\Lambda h}
\end{array}
}. 
\end{align*}
We can now use Lemma \ref{lem:interpolation_inequality}
\begin{equation*}
\begin{aligned}
\av{\widetilde{\sA}_{1}^2}_{s, \lambda} \leqslant & \ \ck_s  \bra{
\av{\partial_1 h}_{s, \lambda} \pare{\av{\mathsf{G}\pare{\Lambda h}}_{0, \lambda} + 1} + 
\av{\partial_1 h}_{0, \lambda} \av{\mathsf{G}\pare{\Lambda h}}_{s, \lambda}
} , \\
%-------------------------------------------------
\av{\widetilde{	\sA}_2^2}_{s, \lambda} \leqslant & \  \av{\mathsf{G}\pare{\Lambda h}}_{s, \lambda},  
\end{aligned}
\end{equation*}
and control the contributions provided by $ \mathsf{G} $ with Lemma \ref{lem:composition_Wiener_FG}
\begin{equation*}
\begin{aligned}
\av{\widetilde{\sA}_{1}^2}_{s, \lambda} \leqslant &  \ \ck_s  \bra{
1+  \av{h}_{1, \lambda} \pare{
 \frac{1}{1-\ck_0  \av{h}_{1, \lambda}} + \frac{1}{1-\ck_s  \av{h}_{1, \lambda}}
}
}  \av{\Lambda h}_{s, \lambda}, \\
%--------------------------------------
\av{\widetilde{\sA}_2^2}_{s, \lambda} \leqslant & \ \frac{1}{1-\ck_s \av{h}_{1, \lambda}} \ \av{\Lambda h}_{s, \lambda}.
\end{aligned}
\end{equation*}
Using now the smallness asumption \eqref{eq:smllnss1}, we find that
\begin{align}\label{eq:bound_A-I}
\av{\sA\pare{t}-I}_{s, \lambda}\leqslant C(s)\av{\Lambda h\pare{t}}_{s, \lambda}.
\end{align}

\subsection{Estimates for the trace of the velocity potential}
Let us start by considering the energy contributions arising in \eqref{eq:EEpsi1}. We use the product law of Lemma \ref{lem:interpolation_inequality} and the fact that for any $ s, \lambda\geq 0 $ the space $ A^s_\lambda $ is an algebra in order to obtain that
\begin{align*}
\av{\av{\sA^\intercal\nabla\phi}^2}_{1, \mu t} \leqslant & \ 2 \av{\sA^\intercal \nabla\pare{\phi_1+\phi_2}}_{0, \mu t} \av{\sA^\intercal \nabla\pare{\phi_1+\phi_2}}_{1, \mu t}, \\
%----------------------------------------------
\leqslant & \ 2 (\av{\sA-I}_{0, \mu t}+1)\pare{\av{\nabla\phi_1}_{0, \mu t} + \av{\partial_2 \phi_2}_{0, \mu t}} 
(\av{\sA-I}_{1, \mu t}+1)\pare{\av{\nabla\phi_1}_{1, \mu t} + \av{\partial_2 \phi_2}_{1, \mu t}} 
. 
\end{align*}

We combine the above bound with the estimates \eqref{estPHI}, \eqref{estPHI2} and the estimate \eqref{eq:bound_A-I} in order to deduce the inequality
\begin{equation}\label{eq:NLB1}
\av{\av{\sA^\intercal\nabla\phi}^2}_{1, \mu t}  
%----------------------------------------------
 \leqslant  \ C\pare{\av{h}_{1, \mu t }+ 1}^2\av{\xi}_{1, \mu t} 
\av{\Lambda h}_{1, \mu t}
\bra{
 \av{\Lambda \xi}_{1, \mu t} + \pare{\av{h}_{1, \mu t} \av{\Lambda \xi}_{1, \mu t} + \av{\xi}_{1, \mu t} \av{\Lambda h}_{1, \mu t} }
} .
\end{equation}
Invoking the assumption \eqref{eq:smllnss1}, then we can simplify \eqref{eq:NLB1} into
\begin{equation*}
\av{\av{\sA^\intercal\nabla\phi}^2}_{1, \mu t} \leqslant C \av{\Lambda h }_{1, \mu t}\pare{\av{\Lambda h }_{1, \mu t} + \av{\Lambda \xi }_{1, \mu t}}.
\end{equation*}

The next term we want to bound is 
\begin{equation*}
\alpha\av{\frac{\pare{2 + \Lambda h }\Lambda h}{\pare{1+\Lambda h}^2} \ \Lambda^2 \xi}_{1, \mu t}. 
\end{equation*}
First of all we remark that
\begin{equation*}
\frac{\pare{2 + \Lambda h }\Lambda h}{\pare{1+\Lambda h}^2} \ \Lambda^2 \xi = 
\pare{2 + \Lambda h }\pare{1-\mathsf{G}\pare{\Lambda h}}^2\Lambda h  \ \Lambda^2 \xi, 
\end{equation*}
where the function $ \mathsf{G} $ is defined in \eqref{eq:def_G}. We iterate the product rule \eqref{eq:product_rule_Wiener} and we obtain
\begin{multline}
\label{eq:NLB2}
\av{\pare{2 + \Lambda h }\pare{1-\mathsf{G}\pare{\Lambda h}}^2\Lambda h  \ \Lambda^2 \xi }_{1, \mu t} \\
\begin{aligned}
%------------------------------------------
\leq & \  32 \bigg[
\pare{2+\av{h}_{1, \mu t }} \pare{1+\av{\mathsf{G}\pare{\Lambda h}}_{0, \mu t}}^2 \av{h}_{1, \mu t} \av{\Lambda^2 \xi}_{1, \mu t} 
\\
%------------------------------------------
& \qquad + \pare{2+\av{h}_{1, \mu t }} \pare{1+\av{\mathsf{G}\pare{\Lambda h}}_{0, \mu t}}^2 \av{ \Lambda h}_{1, \mu t} \av{\Lambda^2 \xi}_{0, \mu t} \\
%------------------------------------------
& \qquad +  \pare{1+\av{\mathsf{G}\pare{\Lambda h}}_{0, \mu t}}^2 \av{  h}_{1, \mu t} \av{ \Lambda h}_{1, \mu t} \av{\Lambda^2 \xi}_{0, \mu t}\\
%------------------------------------------
& \qquad + 2 \pare{2+\av{h}_{1, \mu t }} \pare{1+\av{\mathsf{G}\pare{\Lambda h}}_{0, \mu t}} \av{\mathsf{G}\pare{\Lambda h}}_{1, \mu t} \av{  h}_{1, \mu t} \av{\Lambda^2 \xi}_{0, \mu t} \bigg]. 
\end{aligned}
\end{multline}
We invoke \ref{lem:composition_Wiener_FG} in order to produce the bound
\begin{align}\label{eq:NLB2.1}
\av{\mathsf{G}\pare{\Lambda h}}_{s, \mu t} \leq \frac{\av{\Lambda h}_{s, \mu t}}{1-\av{h}_{1, \mu  t}},\quad s=0, 1. 
\end{align}
Combining the above estimate with the smallness hypothesis \eqref{eq:smllnss1} for $ h $, inequality \eqref{eq:NLB2} becomes
\begin{equation}\label{eq:NLB3}
\av{\pare{2 + \Lambda h }\pare{1-\mathsf{G}\pare{\Lambda h}}^2\Lambda h  \ \Lambda^2 \xi }_{1, \mu t} \leq  C\av{ \Lambda h}_{1, \mu t}\pare{  \av{\Lambda^2 \xi}_{1, \mu t} + \av{\Lambda^2 \xi}_{0, \mu t} }.
\end{equation}
Thus, we conclude that
\begin{equation}
\label{eq:EEpsi3}
\alpha\av{\frac{\pare{2 + \Lambda h }\Lambda h}{\pare{1+\Lambda h}^2} \ \Lambda^2 \xi}_{1, \mu t} \leq \alpha C\av{ \Lambda h}_{1, \mu t}\av{\Lambda^2 \xi}_{1, \mu t}. 
\end{equation}

Using \ref{lem:interpolation_inequality}, \ref{lem:composition_Wiener_FG} and \ref{lem:product_rule_Wiener_strip}, \eqref{estPHI}, \eqref{estPHI2}, \eqref{eq:elliptic_est_phi1}, \eqref{eq:smllnss1} and \eqref{eq:NLB2.1} we estimate
\begin{align}
\alpha\av{\frac{\phi,_2}{1+\Lambda h} h,_{11}}_{1,\mu t}&\leq \alpha\av{\phi,_2\pare{1-\mathsf{G}\pare{\Lambda h}}}_{0,\mu t}\av{ h,_{11}}_{1,\mu t}+\alpha\av{\phi,_2\pare{1-\mathsf{G}\pare{\Lambda h}}}_{1,\mu t}\av{ h,_{11}}_{0,\mu t}\nonumber\\
&\leq C\alpha(\av{h\pare{t}}_{1, \mu t} + \av{\xi\pare{t}}_{1, \mu t})\bra{ \av{ \Lambda^2 h (t)}_{1, \mu t} +   \av{\Lambda^{2} \xi (t)}_{1, \mu t}  }.\label{estaux2}
\end{align}
The term
$$
\av{\frac{\phi,_2}{1+\Lambda h}\sA^k_j\phi,_k\tilde{n}_j}_{1,\mu t}.
$$
is a lower order contribution and thus, it can be estimated as before. Using interpolation, we find that
\begin{align}
\av{\frac{\phi,_2}{1+\Lambda h}\sA^k_j\phi,_k\tilde{n}_j}_{1,\mu t}
&\leq  C \pare{\av{\Lambda h }_{1, \mu t} + \av{\Lambda \xi }_{1, \mu t}}^2\nonumber\\
&\leq  C\pare{\av{h }_{1, \mu t} + \av{\xi }_{1, \mu t}} \pare{\av{\Lambda^2 h }_{1, \mu t} + \av{\Lambda^2 \xi }_{1, \mu t}}.\label{estaux3}
\end{align}

We estimate now the term 
\begin{equation*}
\alpha\av{\frac{\Lambda^2 h}{\pare{1+\Lambda h}^3} \Lambda \xi}_{1, \mu t }.
\end{equation*}
As before, we consider the identity
\begin{equation*}
\frac{\Lambda^2 h}{\pare{1+\Lambda h}^3} \ \Lambda \xi = \pare{1-\mathsf{G}\pare{\Lambda h}}^3 \Lambda\xi \  \Lambda^2 h.
\end{equation*}
We apply the product rules of Lemma \ref{lem:interpolation_inequality} and we obtain that
\begin{multline*}
\av{ \pare{1-\mathsf{G}\pare{\Lambda h}}^3 \Lambda\xi \  \Lambda^2 h}_{1, \mu t} \\
%---------------------------------------------------
\leq C\pare{1+\av{\mathsf{G}\pare{\Lambda h}}_{0, \mu t}}^2  \Big\{ \pare{1+\av{\mathsf{G}\pare{\Lambda h}}_{0, \mu t}} \pare{\av{\xi}_{1, \mu t}\av{\Lambda^2 h}_{1, \mu t} + \av{\Lambda \xi}_{1, \mu t}\av{\Lambda^2 h}_{0, \mu t}}\\
%---------------------------------------------------
+  \av{\mathsf{G}\pare{\Lambda h}}_{1, \mu t} \av{ \xi}_{1 , \mu t}\av{\Lambda^2 h}_{0, \mu t} \Big\}.
\end{multline*}
Recalling \eqref{eq:NLB2.1} we further compute that
\begin{multline}\label{eq:NLB5}
\av{ \pare{1-\mathsf{G}\pare{\Lambda h}}^3 \Lambda\psi \  \Lambda^2 h}_{1, \mu t} \\
%---------------------------------------------------
\leq C \pare{1+\frac{\av{\Lambda h}_{0, \mu t}}{1-\av{h}_{1, \mu  t}}}^2  \left\{ \pare{1+\frac{\av{\Lambda h}_{0, \mu t}}{1-\av{h}_{1, \mu  t}}} \pare{\av{\xi}_{1, \mu t}\av{\Lambda^2 h}_{1, \mu t} + \av{\Lambda \xi}_{1, \mu t}\av{\Lambda^2 h}_{0, \mu t}}
\right. \\
%---------------------------------------------------
+ \left. \frac{   \av{ \xi}_{1 , \mu t} }{1-\av{h}_{1, \mu  t}}\av{\Lambda h}_{1, \mu t}^2 \right\}.
\end{multline}
We apply the estimate \eqref{eq:technical_interpolation_estimate} and Young inequality to find that
\begin{align*}
\av{\Lambda h}_{1, \mu t}^2 & \leq 8 \av{h}_{1, \mu t}\av{\Lambda^2 h}_{1, \mu t}, \\
%------------------------------------------------------
\av{\Lambda \xi}_{1, \mu t}\av{\Lambda h}_{1 , \mu t}  & \leq 2 \pare{ \av{\xi}_{1, \mu t} + \av{h}_{1 , \mu t} } \pare{ \av{\Lambda^2 \xi}_{1, \mu t} + \av{\Lambda^2 h}_{1 , \mu t} }.
\end{align*}
These inequalities applied to \eqref{eq:NLB5} together with \eqref{eq:smllnss1} lead us to
\begin{equation}
\label{eq:EEpsi4}
\alpha\av{\frac{\Lambda^2 h}{\pare{1+\Lambda h}^3} \ \Lambda \xi}_{1, \mu t} \leq \alpha C\pare{\av{\xi}_{1, \mu t } + \av{h}_{1, \mu t }} \pare{\av{ \Lambda^2 h}_{1, \mu t} +  \av{ \Lambda^2 \xi}_{1, \mu t}} . 
\end{equation}

The next term in \eqref{eq:EEpsi1} that we have to estimate is
\begin{equation*}
\alpha\av{\frac{\phi_2,_{22}}{\pare{1+\Lambda h}^2}  }_{1, \mu t}. 
\end{equation*}
We observe that
\begin{equation*}
\frac{\phi_2,_{22}}{\pare{1+\Lambda h}^2} = \pare{1-\mathsf{G}\pare{\Lambda h}}^2 \phi_2,_{22} ,
\end{equation*}
so that we can apply repeatedly the product law of \ref{lem:interpolation_inequality} in order to obtain
\begin{multline}\label{eq:NLB8}
\av{\pare{1-\mathsf{G}\pare{\Lambda h}}^2 \phi_2,_{22}}_{1, \mu t} \\
%---------------------------------------------------
\leq 8 \bra{\pare{1+\av{\mathsf{G}\pare{\Lambda h}}_{0, \mu t}}^2 \av{\phi_2,_{22}}_{1, \mu t} + 2 \pare{1+\av{\mathsf{G}\pare{\Lambda h}}_{0, \mu t}}\av{\mathsf{G}\pare{\Lambda h}}_{1, \mu t} \av{\phi_2,_{22}}_{0, \mu t} }.
\end{multline}
We apply the result in \ref{lem:composition_Wiener_FG}, the estimates \eqref{eq:technical_interpolation_estimate}, the smallness hypothesis \eqref{eq:smllnss1}, the bounds \eqref{estPHI}, \eqref{estPHI2} and \eqref{eq:NLB2.1} in order to transform \eqref{eq:NLB8} into 
\begin{align}
\alpha\av{\frac{\phi_2,_{22}}{\pare{1+\Lambda h}^2}  }_{1, \mu t} &= 
\alpha\av{\pare{1-\mathsf{G}\pare{\Lambda h}}^2 \phi_2,_{22}}_{1, \mu t}\nonumber
\\ &\leq \alpha C(\av{h\pare{t}}_{1, \mu t} + \av{\xi\pare{t}}_{1, \mu t})\bra{ \av{ \Lambda^2 h (t)}_{1, \mu t} +   \av{\Lambda^{2} \xi (t)}_{1, \mu t}  }\label{eq:EEpsi5}.
\end{align}

The last term to estimate in \eqref{eq:EEpsi1} is
\begin{equation*}
\alpha\av{\frac{\Lambda^2 h}{\pare{1+\Lambda h}^3}\phi_2,_2}_{1, \mu t}. 
\end{equation*}
As usual, we remark that
\begin{equation*}
\frac{\Lambda^2 h}{\pare{1+\Lambda h}^3} \phi_2,_2 = 
\pare{1-\mathsf{G}\pare{\Lambda h}}^3
\Lambda^2 h \ \phi_2,_2.
\end{equation*}
Thus using Lemma \ref{lem:interpolation_inequality}
\begin{multline*}
\av{\pare{1-\mathsf{G}\pare{\Lambda h}}^3
\Lambda^2 h \ \phi_2,_2}_{1, \mu t} \leq 32 \left[ \pare{1+\av{\mathsf{G}\pare{\Lambda h}}_{0, \mu t}}^3 \av{ \phi_2,_2 }_{0, \mu t} \av{\Lambda^2 h}_{1, \mu t} \right.  \\
%--------------------------------------------------------------
+   \pare{1+\av{\mathsf{G}\pare{\Lambda h}}_{0, \mu t}}^3  \av{\Lambda^2 h}_{0, \mu t} \av{\phi_2,_2 }_{1, \mu t}\\
 + 
\left.3  \pare{1+\av{\mathsf{G}\pare{\Lambda h}}_{0, \mu t}}^2 \av{\mathsf{G}\pare{\Lambda h}}_{1, \mu t} \av{\Lambda^2 h}_{0, \mu t} \av{ \phi_2,_2 }_{0, \mu t}\right]. 
\end{multline*}
We invoke  \eqref{eq:elliptic_est_phi2_trace}, \eqref{eq:NLB2.1} and the smallness hypothesis \eqref{eq:smllnss1}, in an analogous way as what has been done for the previous terms. Then we obtain the bound
\begin{equation}\label{eq:EEpsi6}
\alpha\av{\frac{\Lambda^2 h}{\pare{1+\Lambda h}^3} \partial_2\phi_2}_{1, \mu t}  \leq \alpha C(\av{h\pare{t}}_{1, \mu t} + \av{\xi\pare{t}}_{1, \mu t})\bra{ \av{ \Lambda^2 h (t)}_{1, \mu t} +   \av{\Lambda^{2} \xi (t)}_{1, \mu t}  }.
\end{equation}

We insert the estimates \eqref{eq:NLB1}, \eqref{eq:EEpsi3}, \eqref{estaux2} \eqref{eq:EEpsi4}, \eqref{eq:EEpsi5} and \eqref{eq:EEpsi6} in \eqref{eq:EEpsi1} and use interpolation in Wiener spaces to obtain the energy inequality
\begin{align}
\ddt \av{\xi}_{1, \mu t} + \frac{\alpha}{2} \av{\Lambda^2 \xi}_{1, \mu t} &\leq \av{h}_{1, \mu t} +\alpha C(\av{h\pare{t}}_{1, \mu t} \nonumber\\
&\quad+ \av{\xi\pare{t}}_{1, \mu t})\bigg{[} \av{ \Lambda^2 h (t)}_{1, \mu t} +   \av{\Lambda^{2} \xi (t)}_{1, \mu t}  bigg{]}\nonumber\\
&\quad+C (\av{ h(t) }_{1, \mu t}+\av{ \xi(t) }_{1, \mu t})\pare{\av{\Lambda^2 h(t) }_{1, \mu t} + \av{\Lambda^2 \xi(t) }_{1, \mu t}}. \label{eq:EEpsi7}
\end{align}

\subsection{Estimates for the interface}
We provide now suitable nonlinear bounds for \eqref{eq:EEh1}. We start studying the term 
\begin{equation*}
\av{\sA_2^2 \phi,_2}_{1, \mu t} .
\end{equation*}
We use now the fact that the space $ A^1_{\mu t}, \ t\geq 0 $ is an algebra and we get
\begin{equation*}
\av{\sA_2^2 \phi,_2}_{1, \mu t} \leq 
(\av{\sA-I}_{1, \mu t}+1)\pare{\av{ \phi_1,_2}_{1, \mu t} + \av{ \phi_2,_2}_{1, \mu t}}
. 
\end{equation*}
We apply the estimates \eqref{eq:bound_A-I} and \eqref{estPHI}, \eqref{estPHI2} and \eqref{eq:elliptic_est_phi1} and we obtain that
\begin{equation*}
\av{ \sA_2^2 \phi,_2}_{1, \mu t} \leq C
\av{\Lambda h}_{1, \mu t}
\pare{
\av{ \Lambda \xi}_{1, \mu t} +  \pare{     \av{h}_{1, \mu t} \av{\Lambda \xi }_{1 , \mu t} +  \av{ \Lambda h}_{1 , \mu t} \av{ \xi }_{1, \mu t} }}+\av{\Lambda \xi}_{1, \mu t}.
\end{equation*}
We combine the above estimate with the interpolation inequality \eqref{eq:technical_interpolation_estimate} and the smallness hypothesis \eqref{eq:smllnss1} and we deduce the bound
\begin{equation}
\label{eq:EEh2}
\av{ \sA_2^2 \phi,_2}_{1, \mu t} \leq  C(\av{h}_{1, \mu t}+\av{\xi}_{1, \mu t}) \pare{\av{ \Lambda^2 h}_{1, \mu t} +  \av{ \Lambda^2 \xi}_{1, \mu t}}+\av{\Lambda \xi}_{1, \mu t}. 
\end{equation}

The last term we have to deal with is
\begin{equation*}
\av{\sA_1^k \phi,_k \ \partial_1 h}_{1, \mu t}.
\end{equation*}
Applying Lemma \ref{lem:interpolation_inequality}, we obtain that
\begin{equation*}
\av{\sA_1^k\phi,_k \ \partial_1 h}_{1, \mu t} \leq 2 \bra{
\av{\sA_1^k\phi,_k}_{1, \mu t}
 \av{  \partial_1 h}_{0, \mu t} 
+\av{\sA_1^k\phi,_k}_{0, \mu t}
 \av{  \partial_1 h}_{1, \mu t} }. 
\end{equation*}
We use now the fact that the space $ A^{s}_{\mu t}, \ s, t\geq 0, $ is an algebra and we find that
\begin{multline*}
\av{\sA_1^k\phi,_k \ \partial_1 h}_{1, \mu t} \leq C
(\av{\sA-I}_{1, \mu t}+1)\pare{\av{\nabla \phi_1}_{1, \mu t} + \av{\partial_2 \phi_2}_{1, \mu t}}
 \av{  \partial_1 h}_{0, \mu t} 
\\+C(\av{\sA-I}_{0, \mu t}+1)\pare{\av{\nabla \phi_1}_{0, \mu t} + \av{\partial_2 \phi_2}_{0, \mu t}}
 \av{  \partial_1 h}_{1, \mu t}
 .
\end{multline*}
We use the estimates \eqref{estPHI}, \eqref{estPHI2}, \eqref{eq:elliptic_est_phi1},  \eqref{eq:bound_A-I} together with the smallness hypothesis \eqref{eq:smllnss1} to conclude 
\begin{equation}
\label{eq:EEh3}
\av{\sA_1^k\phi,_k \ \partial_1 h}_{1, \mu t} \leq C(\av{h}_{1, \mu t}+\av{\xi}_{1, \mu t}) \pare{\av{ \Lambda^2 h}_{1, \mu t} +  \av{ \Lambda^2 \xi}_{1, \mu t}}. 
\end{equation}

We use now the nonlinear bounds \eqref{eq:EEh2} and \eqref{eq:EEh3} in \eqref{eq:EEh1} in order to obtain
\begin{equation}
\label{eq:EEh4}
\ddt \av{h}_{1, \mu} + \frac{\alpha}{2} \av{\Lambda^2 h}_{1, \mu t} \leq C(\av{h}_{1, \mu t}+\av{\xi}_{1, \mu t}) \pare{\av{ \Lambda^2 h}_{1, \mu t} +  \av{ \Lambda^2 \xi}_{1, \mu t}}+\av{\Lambda \xi}_{1, \mu t}. 
\end{equation}

\subsection{Closing the estimates}

We combine \eqref{eq:EEpsi7} and \eqref{eq:EEh4} obtaining the following differential inequality
\begin{multline*}
\ddt \pare{\av{\xi}_{1, \mu t}  + \av{h}_{1, \mu t}} + \frac{\alpha}{2}\pare{\av{\Lambda^2 \xi}_{1, \mu t}  + \av{\Lambda^2 h}_{1, \mu t}} \leq \av{h}_{1, \mu t}+\av{\Lambda \xi}_{1, \mu t}\\+C(1+\alpha)(\av{h}_{1, \mu t}+\av{\xi}_{1, \mu t}) \pare{\av{ \Lambda^2 h}_{1, \mu t} +  \av{ \Lambda^2 \xi}_{1, \mu t}}.
\end{multline*}
Now we see that, if 
$$
\frac{\alpha}{2}>1,
$$
and
$$
(\av{h_0}_{1}+\av{\xi_0}_{1})\leq C(\alpha)\ll1 
$$
an application of Poincar\'e and Gronwall inequalities allows us to deduce the estimate
\begin{equation*}
\av{\xi\pare{t}}_{1, \mu t}  + \av{h\pare{t}}_{1, \mu t} 
+ \delta\int _0^t\av{\Lambda^2 \xi \pare{t'}}_{1, \mu t'}  + \av{\Lambda^2 h \pare{t'}}_{1, \mu t'} \dt' \leq  \av{\xi_0}_{1}  + \av{h_0}_{1} , 
\end{equation*}
for any $ t > 0 $ and a sufficiently small $\delta=\delta(\alpha).$ This ensures the global in time existence and decay. Furthermore, as $\mu t>0$ for $t>0$, we conclude the instantaneous gain of analyticity for the solution. Equipped with this regularity the uniqueness follows from a standard contradiction argument.

\section*{Acknowledgments}
The research of S.S. is supported by the Spanish Ministry of Economy and Competitiveness MINECO through the project MTM2017-82184-R funded by (AEI/FEDER, UE) and acronym "DESFLU" and by the European Research Council through the Starting Grant project H2020-EU.1.1.-639227 FLUID-INTERFACE. R. G-B has been funded by the grant MTM2017-89976-P from the Spanish Ministry of Economy and Competitiveness MINECO and the grant PID2019-109348GA-I00 from the Spanish Ministry of Science, Innovation and Universities MICIU.

\appendix

\section{Functional spaces}\label{sec:notation}	
We consider the reference interface 
$$ 
\Gamma=\bT 
$$ 
(understood as the interval $[-\pi,\pi]$ with periodic boundary conditions) and the reference domain
$$
\Omega= \bT \times\pare{ -\infty, 0}.
$$
We write $N=(0,1)$ the outward pointing normal to $\Omega$.
\subsection*{Sobolev spaces}	
The $L^p$-based Sobolev spaces (also known as Sobolev-Slobodeckij spaces), $W^{s,p}(\Omega)$, are defined as
$$
W^{s,p}(\Omega)=\left\{u\in L^p(\Omega) \text{ s.t. } \|u\|_{\dot{W}^{s,p}(\Omega)}^p=\|D^{\lfloor s\rfloor} u\|^p_{L^p}+\sum_{|\alpha|=\lfloor s\rfloor}\int_{\Omega}\int_{\Omega}\frac{|D^{\alpha}u(x)-D^{\alpha}u(y)|^p}{|x-y|^{2+(s-\lfloor s \rfloor)p}}dxdy<\infty\right\}.
$$
We write $H^s(\Omega)=W^{s,2}(\Omega)$. Equivalently, we have that $H^s(\bT)$ is the usual $L^2$-based Sobolev spaces defined using the Fourier series as
$$
H^s(\bT)=\left\{h\in L^2(\bT) \text{ s.t. }|h|_{H^s(\bT)}^2=|h|_s^2=\sum_{n}(1+|n|^s)^2|\hat{h}(n)|^2<\infty\right\}.
$$
We remark that we will write
$$
\|u\|_{s}=\|u\|_{H^{s}(\Omega)}.
$$
During the whole work we will make extensive use of the following version of the Trace Theorem
\begin{lemma}Let $s>0.5$ be fixed. Then the trace operator
$$
\tau:H^{s}(\Omega)\rightarrow H^{s-0.5}(\Gamma) 
$$
is bounded
$$
|u|_{s-0.5}\leq C\|u\|_{s}.
$$
Furthermore, the trace operator is surjective.
\end{lemma}

\subsection*{Wiener spaces}	
	Let us define the Wiener space as
	\begin{equation*}
\mathbb{A}^s_\lambda \pare{\bT}=\left\{v\in L^2(\bT)\text{ s.t. }	\av{v}_{\mathbb{A}^s_\lambda \pare{\bT}}= \av{v}_{s, \lambda}  = \sum _{n\in\bZ} \pare{ 1 + \av{n} }^s e^{\lambda\av{n}} \av{\hat{v}\pare{n}}<\infty\right\}. 
	\end{equation*}
We denote $\mathbb{A}^s_0=\mathbb{A}^s$. We observe that if $\lambda>0$, then this space is formed by analytic functions. Moreover the Wiener spaces $ \mathbb{A}^s_\lambda, s, \lambda \geq 0 $ satisfy the following properties:
\begin{lemma}[\cite{GGS19}]\label{lem:interpolation_inequality}
Let $ 0\leq s_1 \leq s_2 $, $ \lambda \geq 0 $ and let $ v\in \mathbb{A}^{s_2}_\lambda $. Let us define, for any $ \theta \in \bra{0, 1} $, $ s_{\theta} = \theta s_1 + \pare{1-\theta}s_2 $, then
\begin{equation}\label{eq:interpolation_inequality}
\av{v}_{s_\theta, \lambda} \leq \av{v}_{s_1, \lambda}^{\theta} \av{v}_{s_2, \lambda}^{1-\theta}, 
\end{equation}
moreover for any $ 0\leq s_1 < s_2 $ we have
\begin{equation*}
\mathbb{A}^{s_2}_{\lambda} \Subset \mathbb{A}^{s_1}_{\lambda}. 
\end{equation*}
Define
$$
\ck_q =\left\{\begin{array}{ccc}1 & \text{ if }&q\in\bra{0, 1}\\
2^{q} & \text{ if }&q>1\end{array}\right. 
$$
and let $  f,  g \in \mathbb{A}^0_\lambda $ be such that $ \Lambda^r f,  \Lambda^r g \in \mathbb{A}^s_\lambda, \ s, r, \lambda \geq 0 $, then
\begin{equation}\label{eq:product_rule_Wiener}
\begin{aligned}
\av{ \Lambda^r \pare{fg}}_{s, \lambda} & \leq \ck_s \ck_r \pare{ \Big. \av{f}_{0, \lambda} \av{ \Lambda^r g}_{s, \lambda} + \av{ \Lambda^r f}_{s, \lambda} \av{g}_{0, \lambda}} , 
\end{aligned}
\end{equation}
When $ s=0 $ the following improved product rule holds
\begin{equation*}
\av{fg}_{0, \lambda}\leq \av{f}_{0, \lambda}\av{g}_{0, \lambda}. 
\end{equation*}
Furthermore, for any $ r, s, \lambda \geq 0, \ n\in \bN, \ n\geq 2 $ let us define
\begin{align*}
K_{r, s, n} = \left\lbrace
\begin{aligned}
& n & \pare{ r,  s } & \in \bra{0, 1}^2,  \\
%-----------------------------------------
 &  \frac{\ck_r \ck_s \pare{ \pare{ \ck_r \ck_s }^{n-1} -1}}{\ck_r \ck_s - 1} & \pare{r, s}& \in \pare{1, \infty}^2 ,
\end{aligned}
\right. &&
%-------------------------------------------------
K_{s, n} = K_{0, s, n}. 
\end{align*}
Then for any $ v\in \mathbb{A}^s_\lambda $
\begin{equation*}
\av{\Lambda^r \pare{ v^n }}_{s, \lambda} \leq K_{r, s, n} \av{v}_{0, \lambda}^{n-1} \av{\Lambda^r v}_{s, \lambda}. 
\end{equation*}
\end{lemma}
The previous Lemma can also be applied to the composition with an analytic function. For instance, we have that
\begin{lemma}[\cite{GGS19}] \label{lem:composition_Wiener_FG}
Let us define the function
\begin{align}\label{eq:def_G}
\mathsf{G}\pare{x} = \frac{x}{1+x}, 
\end{align}
Given $ v \in A^s_{\lambda}, \ s, \lambda \geq 0 $ such that $ \av{v}_{0, \lambda} < \min \set{ 1  ,  \ck_s^{-1}} $, then
\begin{align*}
\av{\mathsf{G} \circ v}_{s, \lambda} \leq &  \frac{1}{1-\ck_s \av{v}_{0, \lambda} } \ \av{v}_{s, \lambda} .
\end{align*}
\end{lemma}	
	
Finally, let us state an interpolation estimate
\begin{lemma}
Let $ f $ be a zero-mean smooth function. Then the following interpolation inequality holds true
\begin{equation}\label{eq:technical_interpolation_estimate}
\av{f }_{s, \lambda}\leq 2^{1+\frac{s}{s+1}} \av{f}_{0, \lambda}^{\frac{1}{s+1}} \av{\Lambda f}_{s, \lambda}^{1-\frac{1}{s+1}}. 
\end{equation}
\end{lemma}
\begin{proof}
The proof is a rather standard interpolation argument that we omit for the sake of brevity.
\end{proof}

\subsection*{Wiener-Sobolev spaces}
	We need to define functional spaces in the horizontally periodic open half-plane $ \Omega $. Let $ s, \lambda \geq 0, \ p\in\bra{1, \infty} $ and $ k\in\bN $. We define the Wiener-Sobolev space 
		\begin{equation}\label{eq:norms_cA}
		\cA^{s, k}_\lambda\pare{\Omega}=\left\{u \text{ s.t. }\norm{u}_{\cA^{s, k}_\lambda} = \sum_n \pare{ 1+ \av{n} }^s e^{\lambda\av{n}} \int_{-\infty}^0 \av{\partial_2^k \hat{u}\pare{n, x_2}} \dx_2 <\infty.\right\}
		\end{equation}
These Banach spaces were first introduced in \cite{GGS19} to study a related free boundary problem focusing on the motion of an interface in a porous medium. We denote $	\cA^{s, k} = \cA^{s, k}_0 .$
	
For the sake of readability, some properties of these spaces are collected in the following Lemma:
\begin{lemma}[\cite{GGS19}]\label{lem:product_rule_Wiener_strip}
The map $ v\mapsto v_{x_2 =0} $ is continuous from $ \cA^{s, 1}_\lambda $ onto $ A^s_\lambda  $ and the following estimate holds
\begin{equation*}
\av{\left. v\right|_{x_2=0}}_{s, \lambda} \leq  \norm{v}_{\cA^{s, 1}_\lambda}.
\end{equation*}
In addition, if $ f, g \in \cA^{s, 1}_\lambda, \ s > 0, \ r, \lambda\geq 0 $, the following product rule holds true
\begin{equation}
\label{eq:product_rule_Wiener_strip}
\norm{\Lambda^r \pare{ fg }}_{\cA^{s, 1}_\lambda} \leq 2 \ck_r \ck_s \pare{\norm{\Lambda^r f}_{\cA^{s, 1}_\lambda}\norm{g}_{\cA^{0, 1}_\lambda} + \norm{f}_{\cA^{0, 1}_\lambda}\norm{\Lambda^r g}_{\cA^{s, 1}_\lambda}},  
\end{equation}
and
\begin{equation*}
\norm{\Lambda^r\pare{f^n}}_{\cA^{s, 1}_\lambda}\leq 2^n K_{r, s, n} \norm{f}^{n-1}_{\cA^{0, 1}_\lambda} \norm{\Lambda^r f}_{\cA^{s, 1}_\lambda}, 
\end{equation*}
where $K_{r,s,n}$ was defined in Lemma \ref{lem:interpolation_inequality}.
If $ r=s=0 $
\begin{equation*}
\norm{fg}_{\cA^{0, 1}_\lambda} \leq 2 \norm{f}_{\cA^{0, 1}_\lambda}\norm{g}_{\cA^{0, 1}_\lambda} .  
\end{equation*}
Furthermore, if $ f, g \in \cA^{s, 2}_\lambda $, then
\begin{equation*}
\norm{\Lambda^r \pare{ fg }}_{\cA^{s, 2}_{\lambda}}\leq 2 \ck_r \ck_s \pare{
\norm{\Lambda^r f}_{\cA^{s, 2}_{\lambda}} \norm{g}_{\cA^{0, 1}_{\lambda}} + \norm{f}_{\cA^{0, 1}_{\lambda}} \norm{\Lambda^r g}_{\cA^{s, 2}_{\lambda}}
+  \norm{f}_{\cA^{0, 2}_{\lambda}} \norm{\Lambda^r g}_{\cA^{s, 1}_{\lambda}} +\norm{\Lambda^r f}_{\cA^{s, 1}_{\lambda}} \norm{g}_{\cA^{0, 2}_{\lambda}}
}. 
\end{equation*}
\end{lemma}

Finally we provide optimal elliptic estimates for pseudo-differential operators in the strip $ \Omega $:
\begin{lemma}[\cite{GGS19}] \label{lem:elliptic_estimate_via_symbol}
Let $ s, \lambda \geq 0   \  f: \bR \to \bR $ be $ W^{j,  1 }_{\loc}\pare{\bR}, \ j\in \bN \setminus \set{0} $ and $ u\in A^{s+j-1}_\lambda $. Let $ F \in W^{1, 1}_{\loc}\pare{\bR} $ be the primitive of $ \av{f^{\pare{j}}} $ 
 and suppose there exists  $ C_f $ a strictly positive constant depending on $ f $ only such that
\begin{align*}
F\pare{0} - F\pare{- \av{n}} \leq C_f, 
\end{align*}
uniformly in  $ n\in\bZ $. 
 Then  one has
\begin{equation*}
\norm{ f\pare{ x_2   \Lambda } u }_{\cA^{s, j}_\lambda} \leq C_f\av{ \Lambda^{j-1} u}_{s, \lambda} .
\end{equation*}
\end{lemma}

We observe that constant $ C_f $ is independent of $ n $.

\section{Elliptic estimates in Wiener-Sobolev spaces}\label{sec:elliptic_estimates_appendix}
In the present section we prove elliptic estimates in Wiener-Sobolev spaces for the following Poisson problem for a given regular $ {\bf g} = \pare{g_1, g_2}^{\intercal} $: 
\begin{equation}
\label{eq:Poisson}
\left\lbrace
\begin{aligned}
& \Delta \varphi = \nabla \cdot {\bf g}, & \text{ in } &\cS,  \\
& \varphi \pare{x_1, 0} =0,  & \text{ on } & \Gamma \\
& \partial_2 \varphi\pare{x_1, x_2} =0, & \text{ as } & x_2 \to -\infty . 
\end{aligned}
\right.
\end{equation}
This result is similar to the one in \cite{GGS19} (see also \cite{GS19}). This proposition is used in the nonlinear energy estimates below.

\begin{prop} \label{prop:elliptic_estimates}
Let $ {\bf g}\in \cA^{s_0, 1}_\lambda , \ s_0\geq 1, \ \lambda\geq 0  $ and let $\varphi$ be the solution to \eqref{eq:Poisson}. Then $\nabla\varphi$ can be recovered from $\textbf{g}$ as integral operator and, furthermore, it satisfies the estimates
\begin{align*}
\norm{\Lambda^r \nabla \varphi}_{\cA^{s, 1}_\lambda} & \leq 12 \norm{\Lambda^r {\bf g}}_{\cA^{s, 1}_\lambda}, \\
%----------------------------------------------------------
\norm{\nabla\varphi}_{\cA^{s, 2}_{\lambda}} & \leq 12 \norm{ \Lambda {\bf g} }_{\cA^{s, 1}_{\lambda}}  + 4 \norm{  {\bf g} }_{\cA^{s, 2}_{\lambda}}
\end{align*}
 for any $ r, s\geq 0 $.
\end{prop}

\begin{proof}
For the sake of readability we divide the proof in several steps. 

\textbf{Step 1: derivation of the explicit solution of \eqref{eq:Poisson} }
We write
\begin{equation*}
b = \nabla \cdot {\bf g} ,
\end{equation*}
We have that $ b\in \cA^{s_0-1, 0} $, and since $ s_0-1 \geq 0 $ for any fixed $ \bar{x}_2\in \pare{-1, 0} $, we can define its (horizontal) Fourier transform 
$$ 
\pare{\hat{b}\pare{k, \bar{x}_2}}_{k\in\bZ}.
$$
The Fourier transform in $ x_1 $ transforms \eqref{eq:Poisson} in the sequence of ODE
\begin{align*}
- k^2 \ \hat{\varphi}\pare{k, x_2 } + \partial_2^2 \hat{\varphi}\pare{k, x_2 } = \hat{b}\pare{k, x_2 } , && k \in \bZ, 
\end{align*}
which is explicitly solvable, and its general solution is given by
\begin{equation*}
\hat{\varphi}\pare{k, x_2} = C_1 \pare{k} e^{\av{k}x_2} + C_2 \pare{k} e^{-  \av{k}x_2} - \int_0^{x_2} \frac{\hat{b}\pare{k, y_2}}{2 \av{k}} \bra{e^{\av{k}\pare{y_2 - x_2 }} - e^{-\av{k}\pare{y_2 - x_2 }}} \textnormal{d} y_2. 
\end{equation*}
Using now the first boundary condition we deduce that
\begin{equation*}
C_2 = -C_1. 
\end{equation*}
Next we compute
\begin{equation*}
\partial_2 \hat{\varphi}\pare{k, x_2} = C_1 \pare{k} \av{k} \pare{ e^{\av{k}x_2} +  e^{- \av{k}x_2}} \\
%------------------------------
 + \frac{1}{2}\int_0^{x_2} {\hat{b}\pare{k, y_2}} \bra{e^{\av{k}\pare{y_2 - x_2 }} + e^{\av{k}\pare{x_2 - y_2 }}} \textnormal{d} y_2 , 
\end{equation*} 
which combined with the condition $ \partial_2 \varphi \pare{x_1, -\infty } =0  $ imply that
\begin{equation*}
C_1 \pare{k} = \frac{1}{2  \av{k} }\int _{-\infty}^0 {\hat{b}\pare{k, y_2}} e^{\av{k}y_2} \textnormal{d} y_2 . 
\end{equation*}
From the previous considerations we find that the unique solution is given by
\begin{equation}
\label{eq:psi_explicit}
\hat{\varphi} \pare{k, x_2} = \frac{1}{2  \av{k} }\int _{-\infty}^0 {\hat{b}\pare{k, y_2}} e^{\av{k}y_2} \textnormal{d} y_2 \  \pare{ e^{\av{k}x_2} -  e^{- \av{k}x_2}} \\
%----------------------------------------------------------------
+ \int_{x_2}^0 \frac{\hat{b}\pare{k, y_2}}{2 \av{k}} \bra{e^{\av{k}\pare{y_2 - x_2 }} - e^{-\av{k}\pare{y_2 - x_2 }}} \textnormal{d} y_2 .
\end{equation}

Let us remark that if we define
\begin{equation}
\label{eq:Pi_j}
\begin{aligned}
\Pi_1 \pare{\av{k}, y_2, x_2} = & \ e^{\av{k}y_2} \ \sinh\pare{\av{k}x_2}, & y_2 \in & \ \pare{-\infty, x_2},  \\
%-----------------------------------------------------------
\Pi_2 \pare{\av{k}, y_2, x_2} = & \ e^{\av{k}y_2} \ \sinh\pare{\av{k}x_2} + \sinh \pare{ \av{k}  \pare{y_2 - x_2} }, & y_2 \in & \ \bra{ x_2, 0},   
\end{aligned}
\end{equation}
then $ \hat{\varphi} $ in \eqref{eq:psi_explicit} can be rewritten as
\begin{equation}
\label{eq:psi_explicit2}
\hat{\varphi} \pare{k, x_2} =  \frac{1}{ \av{k} } \bra{  \int _{-\infty}^{x_2} \Pi_1\pare{\av{k}, y_2, x_2} \hat{b}\pare{k, y_2} \textnormal{d} y_2   \\
%-----------------------------------------------------------
+\int _{x_2}^0  \Pi_2\pare{ \av{k}, y_2, x_2} \hat{b}\pare{k, y_2} \textnormal{d} y_2} . 
%--------------------------------------------
\end{equation}

We can now use the explicit form of the forcing $ b = \nabla \cdot {\bf g} $ to obtain
\begin{multline*}
\hat{\varphi}\pare{k, x_2 } = \frac{1}{\av{k}} \int_{-\infty}^{x_2} \Pi_1 \pare{\av{k}, y_2, x_2} \pare{i k \hat{g}_1 \pare{k, y_2} + \partial_{y_2}\hat{g}_2 \pare{k, y_2}} \dy_2 \\
%-----------------------------------------------------------
+ \frac{1}{\av{k}} \int_{x_2}^{0} \Pi_2 \pare{ \av{k}, y_2, x_2} \pare{i k \hat{g}_1 \pare{k, y_2} + \partial_{y_2}\hat{g}_2 \pare{k, y_2}} \dy_2.
\end{multline*}
We decompose
\begin{equation*}
\varphi = \varphi_1 + \varphi_2,
\end{equation*}
with
\begin{align}
\hat{\varphi_1}\pare{k, x_2 } = & \  
i\sgn \pare{k} \bra{ \int_{-\infty}^{x_2} \Pi_1 \pare{ \av{k}, y_2, x_2} \hat{g}_1 \pare{k, y_2}  \dy_2 +  \int_{x_2}^{0} \Pi_2 \pare{ \av{k}, y_2, x_2} \hat{g}_1 \pare{k, y_2}  \dy_2} , \label{eq:psi1} \\
%-----------------------------------------------
\hat{\varphi}_2\pare{k, x_2 } = & \ \frac{1}{\av{k}} \bra{ \int_{-\infty}^{x_2} \Pi_1 \pare{ \av{k}, y_2, x_2}  \partial_{y_2}\hat{g}_2 \pare{k, y_2} \dy_2 +\int_{x_2}^{0} \Pi_2 \pare{ \av{k}, y_2, x_2}  \partial_{y_2}\hat{g}_2 \pare{k, y_2} \dy_2}. \label{eq:psi1b}  
\end{align}

\textbf{Step 2: Computation of $ \nabla \varphi $ and $ \nabla \varphi,_2  $  }

Our first goal is to write $\nabla \varphi$ as a function of $ {\bf g}$. We use the explicit definition of $ \varphi_1 $ and $ \varphi_2 $, computed in the previous step and provided respectively in \eqref{eq:psi1} and \eqref{eq:psi1b} together with integration by parts in $y_2$ and the explicit form of $\Pi_j$ \eqref{eq:Pi_j} in order to obtain
\begin{equation}
\label{eq:pa2_psi}
\begin{aligned}
\partial_2 \hat{\varphi}_1  \pare{k, x_2} = &  \  i\sgn \pare{k}  \left[ \int_{-\infty}^{x_2} \partial_{x_2}\Pi_1 \pare{ \av{k}, y_2, x_2} \hat{g}_1 \pare{k, y_2}  \dy_2 
 +  \int_{x_2}^{0} \partial_{x_2} \Pi_2 \pare{ \av{k}, y_2, x_2} \hat{g}_1 \pare{k, y_2}  \dy_2 \right] , \\
%-----------------------------------------------
%-----------------------------------------------
\partial_2 \hat{\varphi}_2  \pare{k, x_2} = & \ -  \hat{g}_2\pare{k, x_2}\\
%--------------------------------------------------------
  & \  + \frac{1}{\av{k}} \int _{-\infty}^{x_2} {\hat{g}_2\pare{k, y_2}}  \partial_{x_2}\partial_{y_2} \Pi_1\pare{\av{k}, y_2, x_2} \textnormal{d} y_2 + \frac{1}{\av{k}}  \int _{x_2}^0 {\hat{g}_2 \pare{k, y_2}} \partial_{x_2}\partial_{y_2} \Pi_2\pare{\av{k}, y_2, x_2} \textnormal{d} y_2 .
\end{aligned}
\end{equation}
Similarly, we deduce 
\begin{equation}
\label{eq:pa1_psi}
\begin{aligned}
 \widehat{\partial_1 \varphi_1} \pare{k, x_2}  = & \  
 i \av{k} \bra{ \int_{-\infty}^{x_2} \Pi_1 \pare{ \av{k}, y_2, x_2} \hat{g}_1 \pare{k, y_2}  \dy_2 +  \int_{x_2}^{0} \Pi_2 \pare{ \av{k}, y_2, x_2} \hat{g}_1 \pare{k, y_2}  \dy_2} ,  \\
%-----------------------------------------------
\widehat{\partial_1 \varphi_2}\pare{k, x_2 } = & \ i\sgn\pare{k} \bra{  \int _{-\infty}^{x_2}   \partial_{y_2} \Pi_1\pare{\av{k}, y_2, x_2}{\hat{g}_2\pare{k, y_2}} \textnormal{d} y_2 
 %--------------------------------------------------------
 + \int _{x_2}^0 \partial_{y_2} \Pi_2\pare{ \av{k}, y_2, x_2} \hat{g}_2 \pare{k, y_2} \textnormal{d} y_2} .
\end{aligned}
\end{equation}

Using the explicit form of $\Pi_j$ \eqref{eq:Pi_j}, we can write 
\begin{equation}
\begin{small}
\label{eq:nd_psi}
\begin{aligned}
\partial_2 \hat{\varphi}_1  \pare{k, x_2} = &  \  \frac{i\sgn \pare{k} }{ \av{k}^2} \left[ \int_{-\infty}^{x_2}    \partial_{y_2}^2 \partial_{x_2} \Pi_1 \pare{ \av{k}, y_2, x_2} \hat{g}_1 \pare{k, y_2}  \dy_2 +  \int_{x_2}^{0}   \partial_{y_2}^2 \partial_{x_2} \Pi_2 \pare{ \av{k}, y_2, x_2} \hat{g}_1 \pare{k, y_2}  \dy_2 \right]  , \\
%-----------------------------------------------
%-----------------------------------------------
\partial_2 \hat{\varphi}_2  \pare{k, x_2} = & \ -  \hat{g}_2\pare{k, x_2}\\
%--------------------------------------------------------
  & \  + \frac{1}{\av{k}} \int _{-\infty}^{x_2} {\hat{g}_2\pare{k, y_2}}  \partial_{x_2}\partial_{y_2} \Pi_1\pare{\av{k}, y_2, x_2} \textnormal{d} y_2 + \frac{1}{\av{k}}  \int _{x_2}^0 {\hat{g}_2 \pare{k, y_2}} \partial_{x_2}\partial_{y_2} \Pi_2\pare{\av{k}, y_2, x_2} \textnormal{d} y_2  , 
  \end{aligned}
  \end{small}
\end{equation}
\begin{equation}
\begin{small}
\label{eq:nd_ps}
\begin{aligned}
  %-----------------------------------------------
\widehat{\partial_1 \varphi_1} \pare{k, x_2}  = & \  \frac{i}{ \av{k}}
 \left[ \int_{-\infty}^{x_2}   \partial_{y_2}^2  \Pi_1 \pare{ \av{k}, y_2, x_2} \hat{g}_1 \pare{k, y_2}  \dy_2  +  \int_{x_2}^{0}  \partial_{y_2}^2  \Pi_2 \pare{ \av{k}, y_2, x_2} \hat{g}_1 \pare{k, y_2}  \dy_2\right] ,  \\
%-----------------------------------------------
 \widehat{\partial_1 \varphi_2}\pare{k, x_2 } = & \  \ i\sgn\pare{k} \bra{  \int _{-\infty}^{x_2}   \partial_{y_2} \Pi_1\pare{\av{k}, y_2, x_2}{\hat{g}_2\pare{k, y_2}} \textnormal{d} y_2 
 %--------------------------------------------------------
 + \int _{x_2}^0 \partial_{y_2} \Pi_2\pare{ \av{k}, y_2, x_2} \hat{g}_2 \pare{k, y_2} \textnormal{d} y_2} .
\end{aligned}
\end{small}
\end{equation}
Thus, we have obtained an expression for $\nabla \varphi$ in terms of $ {\bf g}$ instead of derivatives of ${\bf g}$.

Before taking an extra derivative in $x_2$ we have to rearrange the previous expressions. We can integrate by parts in order to commute the operator $ \partial_{y_2} $ onto $ {\bf g} $. This gives an equivalent expression for $\nabla\varphi$. We find that
\begin{equation}
\label{eq:pa2_psi_2}
\begin{aligned}
\partial_2 \hat{\varphi}_1  \pare{k, x_2} = &  \  \frac{i\sgn \pare{k} }{ \av{k}^2} \left[ \av{k}^2 e^{\av{k}x_2 }\hat{g}_1\pare{k, 0} +  \int_{-\infty}^{x_2}    \partial_{y_2} \partial_{x_2} \Pi_1 \pare{ \av{k}, y_2, x_2}  \partial_{y_2}\hat{g}_1 \pare{k, y_2}  \dy_2 \right. \\
%-----------------------------------------------
& \qquad \qquad \left. +  \int_{x_2}^{0}   \partial_{y_2} \partial_{x_2} \Pi_2 \pare{ \av{k}, y_2, x_2}  \partial_{y_2} \hat{g}_1 \pare{k, y_2}  \dy_2 \right]  , \\
%-----------------------------------------------
%-----------------------------------------------
\partial_2 \hat{\varphi}_2  \pare{k, x_2} = & \ - 2 \hat{g}_2\pare{k, x_2}\\
%--------------------------------------------------------
  & \  + \frac{1}{\av{k}} \int _{-\infty}^{x_2} \partial_{x_2} \Pi_1\pare{\av{k}, y_2, x_2} {\partial_{y_2} \hat{g}_2\pare{k, y_2}}   \textnormal{d} y_2\\
%--------------------------------------------------------
  & \   + \frac{1}{\av{k}}  \int _{x_2}^0 \partial_{x_2} \Pi_2\pare{\av{k}, y_2, x_2} {\partial_{y_2} \hat{g}_2 \pare{k, y_2}}  \textnormal{d} y_2  .
   \end{aligned}
   \end{equation}
Similarly, we derive   
   \begin{align}
 \widehat{\partial_1 \varphi_1} \pare{k, x_2}  = & \
 - i\hat{g}_1 \pare{k, x_2 }+i e^{\av{k} x_2}\hat{g}_1 \pare{k, 0 }\nonumber  \\
& \  +
  \frac{i}{ \av{k}}
 \bigg{[} \int_{-\infty}^{x_2}   \partial_{y_2}  \Pi_1 \pare{ \av{k}, y_2, x_2} \partial_{y_2} \hat{g}_1 \pare{k, y_2}  \dy_2  \nonumber\\
 &\quad+  \int_{x_2}^{0}  \partial_{y_2}  \Pi_2 \pare{ \av{k}, y_2, x_2} \partial_{y_2} \hat{g}_1 \pare{k, y_2}  \dy_2\bigg{]} ,    \label{eq:pa1_psi_2}
 \\
%-----------------------------------------------
\ \widehat{\partial_1 \varphi_2}\pare{k, x_2 } = & \ i\sgn\pare{k} \Bigg[  \int _{-\infty}^{x_2} \partial_{y_2}\hat{g}_2\pare{k, y_2}   \Pi_1\pare{ \av{k}, y_2, x_2} \textnormal{d} y_2\nonumber
 \\&\quad+ \int _{x_2}^0 \partial_{y_2}\hat{g}_2 \pare{k, y_2}   \Pi_2\pare{\av{k}, y_2, x_2} \textnormal{d} y_2 \Bigg] .    \label{eq:pa1_psi_2b}
\end{align}

Taking an $x_2$ derivative of \eqref{eq:pa2_psi_2} and \eqref{eq:pa1_psi_2}  and using \eqref{eq:Pi_j}, we find that
\begin{align}
\partial_2^2 \hat{\varphi}_1  \pare{k, x_2} = &  \   \frac{i\sgn \pare{k} }{ \av{k}^2} \left[ 
\av{k}^3 e^{\av{k}x_2}\hat{g}_1\pare{k, 0}
 +  \int_{-\infty}^{x_2}    \partial_{y_2} \partial_{x_2}^2 \Pi_1 \pare{ \av{k}, y_2, x_2}  \partial_{y_2}\hat{g}_1 \pare{k, y_2}  \dy_2 \right.\nonumber \\
%-----------------------------------------------
& \qquad \qquad \left. +  \int_{x_2}^{0}   \partial_{y_2} \partial_{x_2}^2 \Pi_2 \pare{ \av{k}, y_2, x_2}  \partial_{y_2} \hat{g}_1 \pare{k, y_2}  \dy_2 \right], \label{eq:pa22_psi_1}\\
%-----------------------------------------------
%-----------------------------------------------
\partial_2^2 \hat{\varphi}_2  \pare{k, x_2} = & \- 2\partial_2  \hat{g}_2\pare{k, x_2} + \frac{1}{\av{k}} \int _{-\infty}^{x_2} \partial_{x_2}^2 \Pi_1\pare{\av{k}, y_2, x_2} {\partial_{y_2} \hat{g}_2\pare{k, y_2}}   \textnormal{d} y_2\nonumber\\
%--------------------------------------------------------
  & \   + \frac{1}{\av{k}}  \int _{x_2}^0 \partial_{x_2}^2  \Pi_2\pare{\av{k}, y_2, x_2} {\partial_{y_2} \hat{g}_2 \pare{k, y_2}}  \textnormal{d} y_2   \label{eq:pa22_psi_1b},\\ 
 \widehat{\partial_2 \partial_1 \varphi_1} \pare{k, x_2}  = & \  
- 2i\partial_2 \hat{g}_1 \pare{k, x_2 }+i\av{k} e^{\av{k} x_2}\hat{g}_1 \pare{k, 0 }  \nonumber\\
& \  +
  \frac{i}{ \av{k}}
 \bigg{[} \int_{-\infty}^{x_2} \partial_{x_2}  \partial_{y_2}  \Pi_1 \pare{ \av{k}, y_2, x_2} \partial_{y_2} \hat{g}_1 \pare{k, y_2}  \dy_2 \nonumber\\&\quad +  \int_{x_2}^{0}  \partial_{x_2} \partial_{y_2}  \Pi_2 \pare{ \av{k}, y_2, x_2} \partial_{y_2} \hat{g}_1 \pare{k, y_2}  \dy_2\bigg{]}  ,     \label{eq:pa21_psi}
\\
%-----------------------------------------------
 \widehat{\partial_2 \partial_1 \varphi_2}\pare{k, x_2 } = & \  i\sgn\pare{k} \Bigg[  \int _{-\infty}^{x_2} \partial_{x_2}  \Pi_1\pare{ \av{k}, y_2, x_2} \partial_{y_2}\hat{g}_2\pare{k, y_2}  \textnormal{d} y_2\nonumber\\
&\quad + \int _{x_2}^0  \partial_{x_2} \Pi_2\pare{\av{k}, y_2, x_2}  \partial_{y_2}\hat{g}_2 \pare{k, y_2}   \textnormal{d} y_2 \Bigg].   \label{eq:pa21_psib}
\end{align}

\textbf{Step 3: Elliptic estimates} We take the absolute value of \eqref{eq:pa22_psi_1} and \eqref{eq:pa12_psi} and integrate in $x_2$ to find

\begin{equation}
\label{eq:pa22_psi_2}
\begin{small}
\begin{aligned}
\int _{-\infty}^{0}\av{\partial_2^2 \hat{\varphi}_1  \pare{k, x_2}} \dx_2  \leq &  \ \int _{-\infty}^{0}  \left\lbrace  \Bigg.  \frac{1}{\av{k}^2} \  \norm{ \mathbf{1} _{\bra{-\infty, x_2}}\pare{\cdot}\partial_{y_2} \partial_{x_2}^2 \Pi_1 \pare{\av{k}, \cdot , x_2} }_{L^\infty_{y_2}}\norm{\partial_{y_2}\hat{g}_1 \pare{k, \cdot }}_{L^1_{y_2}}  \right.  \\
%-----------------------------------------------
&  \qquad \qquad +   \frac{1}{ \av{k}^2} \  \norm{ \mathbf{1} _{\bra{ x_2, 0}}\pare{\cdot} \partial_{y_2} \partial_{x_2}^2 \Pi_2 \pare{ \av{k}, \cdot , x_2} }_{L^\infty_{y_2}}\norm{\partial_{y_2}\hat{g}_1 \pare{k, \cdot }}_{L^1_{y_2}}    \\
%-----------------------------------------------
& \qquad \qquad  \left.\Bigg.
+  \av{k} \ e^{\av{k}x_2} \av{\hat{g}_1 \pare{k, 0}} 
 \right\rbrace \dx_2 , \\
%-----------------------------------------------
%-----------------------------------------------
\int _{-\infty}^{0} \av{\partial_2^2 \hat{\varphi}_2  \pare{k, x_2}} \dx_2 \leq & \ \int _{-\infty}^{0}  \left\lbrace \Bigg. 2 \av{ \partial_2 \hat{g}_2\pare{k, x_2} }\right. \\
%--------------------------------------------------------
  & \ \qquad\qquad +   \frac{1}{ \av{k}} \norm{ \mathbf{1} _{\bra{-\infty, x_2}}\pare{\cdot}\partial_{x_2}^2 \Pi_1\pare{ \av{k}, \cdot , x_2}}_{L^\infty_{y_2}} \norm{\partial_{y_2}\hat{g}_2\pare{k, \cdot }}_{L^1_{y_2}}  \\
 %--------------------------------------------------------
  & \ \qquad\qquad +  \left. \frac{1}{ \av{k}} \norm{ \mathbf{1} _{\bra{ x_2, 0}}\pare{\cdot} \partial_{x_2}^2 \Pi_2\pare{ \av{k}, \cdot , x_2}}_{L^\infty_{y_2}} \norm{\partial_{y_2}\hat{g}_2 \pare{k, \cdot}}_{L^1_{y_2}}  \right\rbrace \dx_2 .  
   \end{aligned}
\end{small}
   \end{equation}
Thus, we have to find appropriate estimates for the previous terms. From \eqref{eq:Pi_j} we deduce that for any $ \pare{ j, l }\in \bN^2 $
\begin{equation}
\label{eq:Pi_Derivatives}
\begin{aligned} 
 2\ \partial_{x_2}^{j} \partial_{y_2}^{l} \Pi_1 \pare{\av{k}, y_2, x_2 } 
%--------------------------------------------------------
& =\av{k}^{j+l}
\bra{ 
e^{\av{k}\pare{ y_2 + x_2  }} -\pare{-1}^j e^{\av{k}\pare{y_2-x_2}}
} , 
 %--------------------------------------------------------
& y_2 \in \pare{-\infty, x_2}, 
\\
%--------------------------------------------------------
2\ \partial_{x_2}^{j} \partial_{y_2}^{l} \Pi_2 \pare{  \av{k}, y_2, x_2 } 
%--------------------------------------------------------
& = \av{k}^{j+l}
\bra{ 
e^{\av{k}\pare{ y_2 + x_2  }} -\pare{-1}^l e^{-\av{k}\pare{y_2-x_2}}
}, 
%--------------------------------------------------------
& y_2 \in \bra{ x_2, 0}. \\
%--------------------------------------------------------. 
\end{aligned}
\end{equation}
Since we have the constraint $ y_2 \in \pare{-\infty, x_2} $, we need to compute
\begin{equation*}
\int _{y_2 }^{0} e^{\av{k}\pare{y_2 - x_2}} \dx_2 = -\frac{1}{\av{k}} e^{\av{k} y_2} \pare{1- e^{- \av{k} y_2} } \leq \frac{1}{\av{k}}. 
\end{equation*}

We can integrate \eqref{eq:Pi_Derivatives} in $ x_2 $ as done above from which we obtain

\begin{align*}
 \ \int _{y_2}^{0} \av{\partial_{x_2}^{j} \partial_{y_2}^{l} \Pi_1 \pare{\av{k}, y_2, x_2 }} \dx_2  &
%--------------------------------------------------------
 \leq  \av{k}^{j+l-1}, &&
%--------------------------------------------------------
 y_2 \in \pare{-\infty, x_2}, \\
%--------------------------------------------------------
% \\
%--------------------------------------------------------
 \ \int _{-\infty}^{y_2} \av{\partial_{x_2}^{j} \partial_{y_2}^{l} \Pi_2 \pare{ \av{k}, y_2, x_2 }} \dx_2  & \leq   \av{k}^{j+l-1}, &&
%--------------------------------------------------------
 %--------------------------------------------------------
 y_2 \in \pare{x_2, 0}. 
\end{align*}
Thus, we conclude the following bounds valid for any $ \pare{j, l}\in\bN^2 $
   \begin{equation}
   \label{eq:bounds_integral_Pi_j}
   \begin{aligned}
   \int _{-\infty}^{0} \norm{ \mathbf{1} _{\bra{-\infty, x_2}}\pare{\cdot} \partial_{x_2}^j \partial_{y_2}^l \Pi_1 \pare{\av{k}, \cdot , x_2} }_{L^\infty_{y_2}} \dx_2 & \leq  \av{k}^{j+l-1}, \\
   %-----------------------------------------------------
   \int _{-\infty}^{0} \norm{ \mathbf{1} _{\bra{ x_2, 0}}\pare{\cdot}  \partial_{x_2}^j \partial_{y_2}^l \Pi_2 \pare{ \av{k}, \cdot , x_2} }_{L^\infty_{y_2}} \dx_2  & \leq  \av{k}^{j+l-1}. 
   \end{aligned}
   \end{equation}
   
   We can use the bounds \eqref{eq:bounds_integral_Pi_j} together with
   $$
 \av{\hat{g}_1\pare{k, 0}} \leq\int _{-\infty}^{0}\av{\partial_{y_2}\hat{g}_1\pare{k, y_2}}\dd y_2 
   $$ in \eqref{eq:pa22_psi_2} to find that
   \begin{equation}
   \label{eq:BOUND1}   
   \begin{aligned}
   \int _{-\infty}^{0}\av{\partial_2^2 \hat{\varphi}_1  \pare{k, x_2}} \dx_2  & \leq 2 \int _{-\infty}^{0} \av{\partial_2 \hat{g}_1 \pare{k, y_2}} \dy_2 + \av{\hat{g}_1\pare{k, 0}}, \\
   %---------------------------------------------------
    & \leq 3 \int _{-\infty}^{0} \av{\partial_2 \hat{g}_1 \pare{k, y_2}} \dy_2, \\
   %---------------------------------------------------
   \int _{-\infty}^{0}\av{\partial_2^2 \hat{\varphi}_2  \pare{k, x_2}} \dx_2  & \leq 4 \int _{-\infty}^{0} \av{\partial_2 \hat{g}_2 \pare{k, y_2}} \dy_2 . 
   \end{aligned}
   \end{equation}

We can argue similarly as above in order to obtain that  
\begin{equation}
   \label{eq:pa12_psi}
   \begin{small}
      \begin{aligned}
  \int _{-\infty}^{0} \av{\widehat{\partial_2 \partial_1 \varphi_1} \pare{k, x_2}} \dx_2   \leq & \  
  \int _{-\infty}^{0} \Bigg\lbrace 
 \av{k}e^{\av{k}x_2} \ \av{\hat{g}_1 \pare{k, 0}} +  2 \av{\partial_2 \hat{g}_1\pare{k, x_2}}
\\
%-----------------------------------------------
& \ \qquad\qquad +  \frac{1}{ \av{k}} \  \norm{\mathbf{1} _{\bra{-\infty,  x_2}}\pare{\cdot} \partial_{x_2} \partial_{y_2}  \Pi_1 \pare{ \av{k}, \cdot , x_2}}_{L^\infty_{y_2}}  \norm{\partial_{y_2}\hat{g}_1 \pare{k, \cdot}}_{L^1_{y_2}} \\
%-----------------------------------------------
& \ \qquad\qquad +  \frac{1}{ \av{k}} \  \norm{ \mathbf{1} _{\bra{ x_2, 0}}\pare{\cdot}  \partial_{x_2} \partial_{y_2}  \Pi_2 \pare{ \av{k}, \cdot , x_2}}_{L^\infty_{y_2}}  \norm{\partial_{y_2} \hat{g}_1 \pare{k, \cdot}}_{L^1_{y_2}} \Bigg\rbrace \dx_2 ,  \\
%-----------------------------------------------
  \int _{-\infty}^{0}\av{ \widehat{\partial_2 \partial_1 \varphi_2}\pare{k, x_2 }} \dx_2  \leq & \ \int _{-\infty}^{0} \Bigg\lbrace \norm{ \mathbf{1} _{\bra{-\infty,  x_2}}\pare{\cdot}  \partial_{x_2} \Pi_1\pare{ \av{k}, \cdot , x_2}}_{L^\infty_{y_2}} \norm{\partial_{y_2}\hat{g}_2\pare{k, \cdot }}_{L^1_{y_2}} \\
%-----------------------------------------------
& \qquad\qquad  + \norm{ \mathbf{1} _{\bra{ x_2, 0}}\pare{\cdot}  \partial_{x_2} \Pi_2\pare{ \av{k}, \cdot , x_2}}_{L^\infty_{y_2}} \norm{\partial_{y_2}\hat{g}_2\pare{k, \cdot }}_{L^1_{y_2}} \Bigg\rbrace \dx_2  , 
\end{aligned}
   \end{small}
\end{equation} 
and we use \eqref{eq:bounds_integral_Pi_j} to find the estimate
\begin{equation}
\label{eq:BOUND2}
\begin{aligned}
  \int _{-\infty}^{0} \av{\widehat{\partial_2 \partial_1 \varphi_1} \pare{k, x_2}} \dx_2 & \leq  \av{\hat{g}_1 \pare{k, 0}} + 2 \int _{-\infty}^{0}\av{\partial_2\hat{g}_1 \pare{k, x_2}} \dx_2, \\
%-------------------------------------------------------
& \leq 3 \int _{-\infty}^{0}\av{\partial_2\hat{g}_1 \pare{k, x_2}} \dx_2, \\
%-------------------------------------------------------
  \int _{-\infty}^{0} \av{\widehat{\partial_2 \partial_1 \varphi_2} \pare{k, x_2}} \dx_2 &\leq  2 \int _{-\infty}^{0}\av{\partial_2 \hat{g}_2 \pare{k, x_2}} \dx_2, 
\end{aligned}
\end{equation}
We can now combine the results in \eqref{eq:BOUND1} and \eqref{eq:BOUND2} together with the remark that
\begin{equation*}
\norm{\Lambda^r {\bf g}}_{\cA^{s, 1}_\lambda} = \norm{\Lambda^r g_1}_{\cA^{s, 1}_\lambda} + \norm{\Lambda^r g_2}_{\cA^{s, 1}_\lambda}, 
\end{equation*}
and we finally conclude the desired bound 
\begin{equation*}
\norm{\Lambda^r \nabla \varphi}_{\cA^{s, 1}_\lambda} \leq 12 \norm{\Lambda^r {\bf g}}_{\cA^{s, 1}_\lambda}. 
\end{equation*}

\textbf{Step 5: Higher order estimates} First we have to compute $ \partial_2^2\nabla\varphi $: To do this we take an $ x_2 $ derivative of \eqref{eq:pa22_psi_1} and \eqref{eq:pa21_psi} and use the explicit form of $\Pi_j$ given by \eqref{eq:Pi_j} obtaining 
\begin{equation}
\label{eq:pa222_psi_1}
\begin{small}
\begin{aligned}
\partial_2^3 \hat{\varphi}_1  \pare{k, x_2} = &  \   \frac{i\sgn \pare{k} }{ \av{k}^2} \left[ 
\av{k}^4 e^{\av{k}x_2}\hat{g}_1\pare{k, 0} - \av{k}^3 \partial_2 \hat{g}_1\pare{k, x_2} \bigg. \right. 
\\
 & \qquad \qquad+  \int_{-\infty}^{x_2}    \partial_{y_2} \partial_{x_2}^3 \Pi_1 \pare{ \av{k}, y_2, x_2}  \partial_{y_2}\hat{g}_1 \pare{k, y_2}  \dy_2 \\
%-----------------------------------------------
& \qquad \qquad \left. +  \int_{x_2}^{0}   \partial_{y_2} \partial_{x_2}^3 \Pi_2 \pare{ \av{k}, y_2, x_2}  \partial_{y_2} \hat{g}_1 \pare{k, y_2}  \dy_2 \right], \\
%-----------------------------------------------
%-----------------------------------------------
\partial_2^3 \hat{\varphi}_2  \pare{k, x_2} = & \- 2\partial_2^2  \hat{g}_2\pare{k, x_2}\\
%--------------------------------------------------------
  & \  + \frac{1}{\av{k}} \int _{-\infty}^{x_2} \partial_{x_2}^3 \Pi_1\pare{\av{k}, y_2, x_2} {\partial_{y_2} \hat{g}_2\pare{k, y_2}}   \textnormal{d} y_2\\
%--------------------------------------------------------
  & \   + \frac{1}{\av{k}}  \int _{x_2}^0 \partial_{x_2}^3  \Pi_2\pare{\av{k}, y_2, x_2} {\partial_{y_2} \hat{g}_2 \pare{k, y_2}}  \textnormal{d} y_2   , 
   \end{aligned}
\end{small}
   \end{equation}
   and
\begin{equation}
   \label{eq:pa221_psi}
   \begin{small}
   \begin{aligned}
 \widehat{\partial_2^2  \partial_1 \varphi_1} \pare{k, x_2}  = & \  
\  
- 2i\partial_2^2 \hat{g}_1 \pare{k, x_2 }+i\av{k}^2 e^{\av{k} x_2}\hat{g}_1 \pare{k, 0 }  \\
& \  +
  \frac{i}{ \av{k}}
 \bigg{[} \int_{-\infty}^{x_2} \partial_{x_2}^2  \partial_{y_2}  \Pi_1 \pare{ \av{k}, y_2, x_2} \partial_{y_2} \hat{g}_1 \pare{k, y_2}  \dy_2  \nonumber\\&\quad+  \int_{x_2}^{0}  \partial_{x_2}^2 \partial_{y_2}  \Pi_2 \pare{ \av{k}, y_2, x_2} \partial_{y_2} \hat{g}_1 \pare{k, y_2}  \dy_2\bigg{]}   ,  \\
%-----------------------------------------------
 \widehat{\partial_2^2 \partial_1 \varphi_2}\pare{k, x_2 } = & ik \partial_2 \hat{g}_2\pare{k, x_2} \\
%-----------------------------------------------
&  \ + i\sgn\pare{k} \Bigg[  \int _{-\infty}^{x_2} \partial_{x_2}^2  \Pi_1\pare{ \av{k}, y_2, x_2} \partial_{y_2}\hat{g}_2\pare{k, y_2}  \textnormal{d} y_2
\nonumber\\&\quad + \int _{x_2}^0  \partial_{x_2}^2 \Pi_2\pare{\av{k}, y_2, x_2}  \partial_{y_2}\hat{g}_2 \pare{k, y_2}   \textnormal{d} y_2 \Bigg] , 
\end{aligned}
   \end{small}
\end{equation} 
Integrating \eqref{eq:pa222_psi_1} and \eqref{eq:pa221_psi} in $x_2$, we deduce that
   \begin{equation}
\label{eq:pa222_psi_2_bound}
\begin{small}
\begin{aligned}
\int _{-\infty}^{0}\av{\partial_2^3 \hat{\varphi}_1  \pare{k, x_2}} \dx_2  \leq &  \ \int _{-\infty}^{0}  \left\lbrace  \Bigg.  \frac{1}{\av{k}^2} \  \norm{ \mathbf{1} _{\bra{-\infty, x_2}}\pare{\cdot}\partial_{y_2} \partial_{x_2}^3 \Pi_1 \pare{\av{k}, \cdot , x_2} }_{L^\infty_{y_2}}\norm{\partial_{y_2}\hat{g}_1 \pare{k, \cdot }}_{L^1_{y_2}}  \right.  \\
%-----------------------------------------------
&  \qquad \qquad +   \frac{1}{ \av{k}^2} \  \norm{ \mathbf{1} _{\bra{ x_2, 0}}\pare{\cdot} \partial_{y_2} \partial_{x_2}^3 \Pi_2 \pare{ \av{k}, \cdot , x_2} }_{L^\infty_{y_2}}\norm{\partial_{y_2}\hat{g}_1 \pare{k, \cdot }}_{L^1_{y_2}}    \\
%-----------------------------------------------
& \qquad \qquad  \left.\Bigg.
+  \av{k}^2 \ e^{\av{k}x_2} \av{\hat{g}_1 \pare{k, 0}}  + \av{k}\av{\hat{g}_1 \pare{k, x_2 }}
 \right\rbrace \dx_2 , \\
%-----------------------------------------------
%-----------------------------------------------
\int _{-\infty}^{0} \av{\partial_2^3 \hat{\varphi}_2  \pare{k, x_2}} \dx_2 \leq & \ \int _{-\infty}^{0}  \left\lbrace \Bigg. 2 \av{ \partial_2^2 \hat{g}_2\pare{k, x_2} }\right. \\
%--------------------------------------------------------
  & \ \qquad\qquad +   \frac{1}{ \av{k}} \norm{ \mathbf{1} _{\bra{-\infty, x_2}}\pare{\cdot}\partial_{x_2}^3 \Pi_1\pare{ \av{k}, \cdot , x_2}}_{L^\infty_{y_2}} \norm{\partial_{y_2}\hat{g}_2\pare{k, \cdot }}_{L^1_{y_2}}  \\
 %--------------------------------------------------------
  & \ \qquad\qquad +  \left. \frac{1}{ \av{k}} \norm{ \mathbf{1} _{\bra{ x_2, 0}}\pare{\cdot} \partial_{x_2}^3 \Pi_2\pare{ \av{k}, \cdot , x_2}}_{L^\infty_{y_2}} \norm{\partial_{y_2}\hat{g}_2 \pare{k, \cdot}}_{L^1_{y_2}}  \right\rbrace \dx_2 .  
   \end{aligned}
\end{small}
   \end{equation}
   \begin{equation}
   \label{eq:pa122_psi_bound}
   \begin{small}
      \begin{aligned}
  \int _{-\infty}^{0} \av{\widehat{\partial_2^2 \partial_1 \varphi_1} \pare{k, x_2}} \dx_2   \leq & \  
  \int _{-\infty}^{0} \Bigg\lbrace 
 \av{k}^2 e^{\av{k}x_2} \ \av{\hat{g}_1 \pare{k, 0}} +  2 \av{\partial_2^2 \hat{g}_1\pare{k, x_2}}
\\
%-----------------------------------------------
& \ \qquad\qquad +  \frac{1}{ \av{k}} \  \norm{\mathbf{1} _{\bra{-\infty,  x_2}}\pare{\cdot} \partial_{x_2}^2 \partial_{y_2}  \Pi_1 \pare{ \av{k}, \cdot , x_2}}_{L^\infty_{y_2}}  \norm{\partial_{y_2}\hat{g}_1 \pare{k, \cdot}}_{L^1_{y_2}} \\
%-----------------------------------------------
& \ \qquad\qquad +  \frac{1}{ \av{k}} \  \norm{ \mathbf{1} _{\bra{ x_2, 0}}\pare{\cdot}  \partial_{x_2}^2 \partial_{y_2}  \Pi_2 \pare{ \av{k}, \cdot , x_2}}_{L^\infty_{y_2}}  \norm{\partial_{y_2} \hat{g}_1 \pare{k, \cdot}}_{L^1_{y_2}} \Bigg\rbrace \dx_2 ,  \\
%-----------------------------------------------
  \int _{-\infty}^{0}\av{ \widehat{\partial_2 \partial_1 \varphi_2}\pare{k, x_2 }} \dx_2  \leq & \ \int _{-\infty}^{0} \Bigg\lbrace \av{k}\av{\partial_2 \hat{g}_2 \pare{k, x_2}} \\
%-----------------------------------------------
& \qquad\qquad\norm{ \mathbf{1} _{\bra{-\infty,  x_2}}\pare{\cdot}  \partial_{x_2}^2 \Pi_1\pare{ \av{k}, \cdot , x_2}}_{L^\infty_{y_2}} \norm{\partial_{y_2}\hat{g}_2\pare{k, \cdot }}_{L^1_{y_2}} \\
%-----------------------------------------------
& \qquad\qquad  + \norm{ \mathbf{1} _{\bra{ x_2, 0}}\pare{\cdot}  \partial_{x_2}^2 \Pi_2\pare{ \av{k}, \cdot , x_2}}_{L^\infty_{y_2}} \norm{\partial_{y_2}\hat{g}_2\pare{k, \cdot }}_{L^1_{y_2}} \Bigg\rbrace \dx_2  . 
\end{aligned}
   \end{small}
\end{equation} 

We apply \eqref{eq:bounds_integral_Pi_j} to get
\begin{equation}
   \label{eq:BOUND3}   
   \begin{aligned}
   \int _{-\infty}^{0}\av{\partial_2^3 \hat{\varphi}_1  \pare{k, x_2}} \dx_2  & \leq  3\av{k} \int _{-\infty}^{0} \av{\partial_2 \hat{g}_1 \pare{k, y_2}} \dy_2 +  \av{k}\av{\hat{g}_1\pare{k, 0}}, \\
   %---------------------------------------------------
    & \leq 4\av{k} \int _{-\infty}^{0} \ \av{\partial_2 \hat{g}_1 \pare{k, y_2}} \dy_2, \\
   %---------------------------------------------------
   \int _{-\infty}^{0}\av{\partial_2^3 \hat{\varphi}_2  \pare{k, x_2}} \dx_2  & \leq 2\av{k} \int _{-\infty}^{0} \av{\partial_2 \hat{g}_2 \pare{k, y_2}} \dy_2 + 2  \int _{-\infty}^{0} \av{\partial_2^2 \hat{g}_2 \pare{k, y_2}} \dy_2 , \\
   %---------------------------------------------------
     \int _{-\infty}^{0} \av{\widehat{\partial_2^2 \partial_1 \varphi_1} \pare{k, x_2}} \dx_2 & \leq  \av{k} \av{\hat{g}_1 \pare{k, 0}}+ 2  \int _{-\infty}^{0} \av{\partial_2^2 \hat{g}_2 \pare{k, y_2}} \dy_2 + 2\av{k} \int _{-\infty}^{0}\av{\partial_2\hat{g}_1 \pare{k, x_2}} \dx_2, \\
%-------------------------------------------------------
& \leq 3\av{k} \int _{-\infty}^{0}\av{\partial_2\hat{g}_1 \pare{k, x_2}} \dx_2  + 2  \int _{-\infty}^{0} \av{\partial_2^2 \hat{g}_2 \pare{k, y_2}} \dy_2  , \\
%-------------------------------------------------------
  \int _{-\infty}^{0} \av{\widehat{\partial_2^2 \partial_1 \varphi_2} \pare{k, x_2}} \dx_2 &\leq  3\av{k} \int _{-\infty}^{0}\av{\partial_2 \hat{g}_2 \pare{k, x_2}} \dx_2. 
   \end{aligned}
   \end{equation}
Collecting the previous bounds we conclude the desired estimate
   \begin{equation*}
   \norm{\nabla\varphi}_{\cA^{s, 2}_{\lambda}} \leq 12 \norm{ \Lambda {\bf g} }_{\cA^{s, 1}_{\lambda}}  + 4 \norm{  {\bf g} }_{\cA^{s, 2}_{\lambda}}. 
   \end{equation*}

\end{proof}

%\section*{References}

\begin{footnotesize}
%\bibliography{references}
%\bibliographystyle{plain}

\end{footnotesize}
\vspace{2cm}

\end{document}